\documentclass[]{article}

\RequirePackage{graphicx, color}
\RequirePackage{caption, subcaption}
\RequirePackage{amscd,amsmath,amsfonts,amssymb}
\RequirePackage{mathrsfs, mathtools}
\RequirePackage{upgreek}
\RequirePackage{fancyhdr}
\RequirePackage{stmaryrd}
\RequirePackage{float}
\RequirePackage{theorem}

\RequirePackage{bm, bbm}
\RequirePackage{a4wide}

\RequirePackage{url}

\graphicspath{%
{./Parts/Images/},%
{./Parts/Images/benchmark/geiger/}
{./Parts/Images/benchmark/geiger/blocking/},%
{./Parts/Images/benchmark/geiger/permeable/},%
{./Parts/Images/benchmark/sotra/}%
{./Parts/Images/transport/}
{./Parts/Images/soultz/}%
}

\newcommand{\jump}[1]{\left\llbracket {#1} \right\rrbracket}

\begin{document}


\def\pathImages{./Parts/Images/}%



\title{Dual Virtual Element Methods\\ for Discrete Fracture Matrix Models}

\author
{
    Alessio
    Fumagalli,
    Eirik
    Keilegavlen\\
    Department of Mathematics, University of Bergen, Bergen, Norway\\
    \url{alessio.fumagalli@uib.no}
    \url{eirik.keilegavlen@uib.no}
}



\maketitle


\abstract
{
        The accurate description of fluid flow and transport in fractured porous
        media is of paramount importance to capture the macroscopic behaviour of
        an oil reservoir, a geothermal system, or a CO$_2$ sequestration site,
        to name few applications.  The construction of accurate simulation model
        for flow in fractures is challenging due to the high ratios between a
        fracture's length and width, which makes modeling by lower-dimensional
        manifolds a natural option.  In this paper we present a
        mixed-dimensional Darcy problem able to describe pressure and Darcy
        velocity in all the dimensions, \textit{i.e.} in the rock matrix, in the
        fractures, and in their intersections. Moreover, we present a
        mixed-dimensional transport problem which, given the Darcy velocity,
        describes coupled advection and diffusion of a passive scalar into the
        fractured porous media.  The approach can handle both conducting and
        blocking fractures.  Our computational grids are created by coarsening
        of simplex tessellations that conform to the fractures surfaces.  An
        accurate choice of the discrete approximation of the previous model, by
        virtual finite element and finite volume, allows us to simulate complex
        problem with a good balance in term of accuracy and computational cost.
        We illustrate the performance of our method by comparing to benchmark
        studies for two-dimensional fractured porous media, as well as a complex
        three-dimensional fracture geometry.
}

\section{Introduction}

Fractures and faults can strongly influence fluid flow in a porous media
acting, depending on their permeability and porosity, as a preferential path or
a barrier. Due to fracture aperture being several orders of magnitude smaller
than any other characteristic sizes in the domain, fracture modeling is one of
the main challenges in subsurface problems.

Geological movements, chemical reaction, or infilling processes may
substantially alter the local orientation and composition of the material
present in the fractures, leading to anisotropies and strong heterogeneities in
both fractures and their intersections. It is thus crucial to be able to include
also these phenomenological aspects in the conceptual model.

Applications where fractures can be determinant for reservoir behaviour include
exploitation of geothermal system, CO$_2$ storage and sequestration, enhance oil
recovery, and nuclear waste disposal. In these fields the solution may
dramatically depends on the presence of fractures, thus a correct derivation of
suitable mathematical models and their accurate numerical solution is essential.

Small scale fractures can be incorporated in the matrix permeability by
analytical or numerical upscaling techniques (\textit{e.g.} see
\cite{Mourzenko2011,Saevik2013,Saevik2014} and
\cite{Karimi-Fard2006,Fumagalli2015,Karimi-Fard2016}), thus only macroscopic
fractures and faults are considered and explicitly described. However, their
number may still be too high for a fully three dimensional representation.
Following the idea presented, among the others, in
\cite{Alboin2000,Martin2005,Amir2005,Tunc2012,Knabner2014,Faille2014a,Schwenck2015,Boon2016}
fractures and faults are represented as lower dimensional objects embedded in
the rock matrix. For this approximation common terminologies are reduced models,
hybrid-dimensional models, and mixed-dimensional models.  Fracture aperture
becomes thus a coefficient in the equations and not a geometrical constrain for
grid generation. The flow equations are averaged along the normal direction of
the fractures obtaining a new set of reduced models suited for the
co-dimensional, that is, lower-dimensional, description of fluid flow.
Suitable coupling conditions are thus necessaries to exchange information
between all the objects (rock matrix, fractures, and their intersections).

The reduction of dimensionality is an essential step to overcome most of the
difficulties associated to the problem. However, in presence of several
fractures which create a complex network more advanced numerical technique are
crucial to obtain an accurate solution in a reasonable amount of time. The main
aspect is the geometrical treatment of the fracture grid with respect to the
rock matrix grid. The common approaches consider conforming and non-conforming
coupling in the case of matching grids or completely a non-matching strategy, in
the latter the rock grid is normally considered as a background and the fracture
grid may adapt to it up to a certain degree.

In the class of conforming methods,
where fracture grids are composed by faces of the rock matrix grid, several
numerical schemes have been considered, ranging from finite element
\cite{Angot2003,Karimi-Fard2003,Martin2005,Angot2009,Chave2017} to finite
volumes \cite{Brenner2015,Brenner2016a,Karimi-Fard2016}, mimetic finite
difference (MFD) \cite{Antonietti,Scotti2017}, and the newly introduced virtual finite
elements (VEM) \cite{Benedetto2014,Fumagalli2016a}. All these methods highlight
specific advantages for example related to local mass conservation, capability
to be implemented in standard software packages, or relax some constraints on
grid cell shapes to name a few.

In the class of non-conforming discretization
the main tool is a mortar coupling between the rock matrix grid and the fracture
grids. The rock mesh is constrained with the position of the fractures but not
strictly with the actual fracture meshes, and vice-versa. Some examples are
reported in \cite{Frih2011,Benedetto2016,Boon2016} where different type of
mortar variables as well as numerical schemes are considered.

Finally, a fully
non matching coupling among the grids requires ad-hoc solutions to establish a
communication between the fractures and the rock matrix. One possibility is
the class of extended finite element methods (XFEM)
\cite{DAngelo2011,Fumagalli2012g,Fumagalli2012d,Berrone2013a,Schwenck2015,DelPra2015a},
where a local enrichment with new basis functions is considered to handle the
non-conformity, or the class of embedded discrete fracture matrix methods (EDFM)
\cite{Li2008,Fumagalli2015,Fumagalli2015d,Tene2016}, where approximate formulae
for fracture-matrix transmissibilities are computed based on geometrical
considerations.

In this paper we consider a conforming discretization with the dual virtual
element approximation to simulate a mixed-dimensional Darcy problem and a finite
volume discretization for the mixed-dimensional transport problem. The first
choice is motivated by the flexibility of the virtual element methods with
respect to the shape of the grid cells. With this choice it is possible to relax
most of the difficulties related with conforming discretization, \textit{e.g.}
the computational cost associated to resolve, by the rock matrix grid,
a complex system of composed by several intersecting fractures. Moreover, this approximation
ensure local mass conservation and is able to consider heterogeneous and
anisotropic permeability tensors. The Darcy velocity is thus suited to be used
in the mixed-dimensional transport problem.  In this case we consider an upwind
scheme, extended to handle the mixed-dimensional nature of the problem.

The paper is structured as follows: in Section \ref{sec:mathematical_model} we
present both the physical equations and the reduced model, with the interface
conditions that couple the matrix-fracture system and the fracture-fracture
system for both the Darcy and transport problems. Section \ref{sec:weak} deals
with the weak and integral formulation of the previous physical processes.
In Section \ref{sec:numerical_discretization} we present the
numerical discretization of the problem with an highlight on the enrichment of
the finite element spaces. In Section \ref{sec:examples} we present
some numerical experiments to asses the effectiveness of the proposed method.
Finally Section \ref{sec:conclusions} is devoted to conclusions and to ongoing
works.




\section{Mathematical model} \label{sec:mathematical_model}

In this section we introduce the mathematical models in the mixed-dimensional
setting, \textit{i.e.} the equations couple different spatial dimensions and are
able to globally describe the quantities of interest. We present two mathematical
models useful for subsurface simulations: the
mixed-dimensional Darcy problem, presented in Subsection \ref{subsec:mix_darcy},
to describe the flow and pressure and the mixed-dimensional transport problem,
presented in Subsection \ref{subsec:mix_transport}, to describe the motion of a
passive scalar transported by the Darcy velocity.

The idea behind the mixed-dimensional formulation is that the fracture aperture
is orders of magnitude smaller than other characteristic sizes of the problem.
Hence, a straightforward mesh construction with the two fracture surfaces as
constraints would produce high numbers of cells and or low-quality cells due to
high aspect ratios or sharp angles.  To avoid this geometrical constraint,
fractures are represented as lower dimensional objects embedded in the rock
matrix.  Fracture intersections, and their intersections again, are considered
as objects of even lower dimensions.  To be specific, in a three-dimensional
domain, we consider fracture surfaces as 2d objects, the intersection between
two fractures form a 1d line, and two intersection lines can meet in a 0d
point.  The physical processes are described via reduced models with suitable
coupling conditions among the objects of different dimensions. The fracture
aperture is now part of the equations and not any more a geometrical constraint.
For more details on the derivation of the reduced model we refer to
\cite{Alboin2000,Martin2005,Amir2005,DAngelo2011,Formaggia2012,Fumagalli2012d,Fumagalli2012g,Berrone2013,Benedetto2014,Faille2014a,Benedetto2016,Fumagalli2016a,Boon2016}.

In the reduced model is a common choice to consider the reduced scalar variables
as averaged and the vector variables as integrated along cross section of the
fracture. In this work we follow this approach.

In the sequel we indicate by $(\cdot, \cdot)_A$ the scalar
product in $L^2(A)$.  Moreover, the trace operator on a domain $A$
will be indicated by $\cdot|_A$.


\subsection{Mixed-dimensional Darcy problem} \label{subsec:mix_darcy}

\begin{figure}
    \centering
    \includegraphics[width=0.25\textwidth]{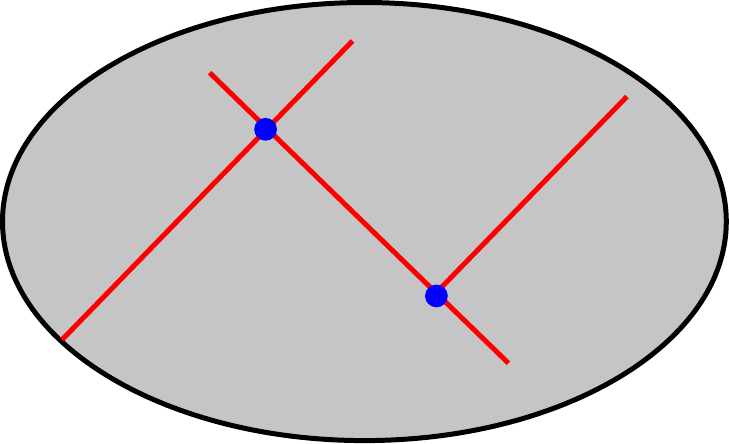}%
    \caption{Domain subdivision for $N=2$. In grey $\Omega^2$, in red
    $\Omega^1$, and in blue $\Omega^0$.}%
    \label{fig:mix_dom}%
\end{figure}

Let us consider a regular domain $\Omega \subset \mathbb{R}^N$, for $N > 0$,
with outer boundary $\partial \Omega$. The domain $\Omega$ is composed by a
single, possibly non-connected, equi-dimensional domain $\Omega^N$ and several
lower-dimensional domains $\Omega^d$ for $d<N$, possibly non-connected. Clearly
$\cup_{d=0, \ldots, N} \Omega^d = \Omega$ and $\Omega^{d} \cap \Omega^{d^\prime}
= \emptyset$ for $d \neq d^{\prime}$. However, we indicate with
$\Gamma^d = \partial \Omega^{d} \cap \Omega^{d-1}$, which is geometrically
equivalent to $\Omega^{d-1}$ but it will be more convenient to keep them separate.
We have $\Gamma = \cup_{d=0, \ldots, N} \Gamma^d$. See
Figure \ref{fig:mix_dom} as an example of the subdivision.

We are interested to the
mixed-dimensional Darcy problem which describes the pressure and velocity fields
in $\Omega$. Following the idea presented in \cite{Boon2016}, we define the
pressure compound $p=(p^0, \ldots, p^N)$ which represents the pressure in each
dimension, indicated as superscript. Similarly we introduce the velocity
compound $\bm{u} = (0, \bm{u}^1, \ldots, \bm{u}^N, 0)$, which is the velocity
field defined in $\Omega^N$ or in the tangential spaces of each dimension $d<N$.
The meaning of the two $0$s will be clarified. In the following problem, each
dimension is coupled with the lower dimension; we consider thus two normals
$\bm{n}$ associated to a $d<N$-dimensional manifold with respect to its ``two
sides'', indicated by $\bm{n}_+$ and $\bm{n}_-$. Given a side, by convention the normal
points from the higher dimensional to the lower dimensional domain.  Moreover,
we indicate by the same subscript $\cdot_\pm$ the restriction of elements in
$\Omega^d$ defined on each side of $\Gamma^{d}$ through a suitable trace
operator.

We consider the Darcy problem for the pressure and velocity as
\begin{subequations} \label{eq:darcy}
\begin{gather}\label{eq:darcy_d}
    \begin{aligned}
        &\bm{u} + K \nabla p = \bm{0}\\
        &\nabla \cdot \check{\bm{u}} - \jump{\hat{\bm{u}} \cdot \hat{\bm{n}}} =
        \check{f}
    \end{aligned}
    \quad \text{in } \Omega.
\end{gather}
The $\hat{\cdot}$ notation means ``on the higher dimensional object'' and the $\check{\cdot}$
means ``on the lower dimensional object'', \textit{i.e.} fix a dimension $d$ the
second of \eqref{eq:darcy_d} becomes
\begin{gather*}
    \nabla \cdot \bm{u}^d - \jump{\bm{u}^{d+1} \cdot \bm{n}^{d+1}} = f^d
    \quad 0 \leq d \leq N.
\end{gather*}
In the system \eqref{eq:darcy_d} the first equation is defined on the same
dimensional domain while the second on mixed-dimensional domain.
In the
previous equation the differential operators are defined accordingly,
\textit{i.e.} for the lower dimensional objects tangential gradient and
divergence. In the special 0d case these objects are the null operator.
$K$ represents the effective (\textit{i.e.} is already
scaled by the ``aperture'' of each dimension) permeability compound $K=(1, K^1,
\ldots, K^N)$, which is a collection of symmetric and positive defined tensors.
For $d<N$, $K^d$ is interpreted as laying in the tangent plane of the relevant physical object.
The source or sink term ${f}$ is defined as the compound
$f=(f^0, \ldots, f^N)$, which are scalar functions defined in each dimension.
Finally, the jump operator $\jump{\hat{\bm{u}} \cdot
\hat{\bm{n}}} = \sum_\pm \hat{\bm{u}} \cdot \hat{\bm{n}}|_\pm$ is done on each
side of $\Omega^{d-1}$ viewed as immersed in $\Omega^d$.

The first equation in \eqref{eq:darcy_d} is the generalized Darcy law in each
dimension, note that for the 0d case it is simply an identity.
The second equation of \eqref{eq:darcy_d} represents the conservation of mass.
The divergence term
describes the conservation of the mass in the dimension and the flow exchange
with the lower dimensional domain (see the derivation of the weak formulation in
Section \ref{sec:weak}), while the second represents the
inflow/outflow from the higher dimension. Due to the definition of the velocity
compound, for $d=N$ equation \eqref{eq:darcy_d} the inflow/outflow from the
``higher dimension'' is null.  While for $d=0$ the second equation represents
the continuity of the fluxes for the $1$-dimensional objects involved in the
intersection.

It is important to note that the mixed-dimensional problem requires that the
quantities defined on $\Omega^{d-1}$ are in communication with the quantities
defined on $\Omega^{d}$ and not directly with quantities defined on
$\Omega^{d+1}$.

The coupling conditions between two subsequent dimensions are defined for the
each side as
\begin{gather}\label{eq:darcy_cc}
    \hat{\bm{u}} \cdot \hat{\bm{n}}|_\pm +\hat{k}_* (\check{p} - \hat{p}|_\pm) = 0
    \quad \text{on } \Gamma.
\end{gather}
The $k_*$ is the effective normal permeability, defined as $k_*=(k_*^0, \ldots,
k_*^{N-1})$.
Clearly, equation \eqref{eq:darcy_cc} is
valid between dimension $d$ and $d+1$ for $0 \leq d < N$. Note that in equation
\eqref{eq:darcy_cc} it is possible to define different effective normal
permeability on each side $\pm$, however, to ease the notation, we assume only one
single value.

Boundary conditions need to be assigned to obtain a well-posed problem. We
indicate by $\partial \Omega^d_{out}$ the (possibly empty) portion of the
boundary of $\Omega^d$ that intersect the boundary of $\Omega$, for $d>0$. We
introduce also $\partial \Omega^d_{in}$ the portion of boundary of $\Omega^d$
that do not intersect the boundary of $\Omega$, for $d>0$. For simplicity we
assume only pressure boundary conditions on $\partial\Omega^d_{out}$, given as
\begin{gather} \label{eq:darcy_bd}
    p^d = \overline{p}^d \quad \text{on } \partial \Omega^d_{out},
\end{gather}
with $\overline{p}^d$ is a given pressure at the boundary for each $d>0$.
On the internal portion $\partial \Omega^d_{in}$, for $d>0$, we assign the
so-called ``tip-condition'', namely
\begin{gather}\label{eq:darcy_int}
    \bm{u}^d \cdot \bm{n}|_{\partial \Omega^d_{in}} = 0 \quad \text{on } \partial
    \Omega^d_{in},
\end{gather}
where in the previous equation $\bm{n}$ stands for the outward unit normal of $\partial
\Omega^d_{in}$.
\end{subequations}

Equation \eqref{eq:darcy} is the mixed-dimensional Darcy problem, formulated in
terms of pressure and velocity.  The problem \eqref{eq:darcy} may be recast into
a pure pressure formulation, however the numerical scheme introduced in Section
\ref{sec:numerical_discretization} considers explicitly the pressure and
velocity fields as unknowns.

Following the idea presented in \cite{Frih2011,Boon2016,Boon2017}, it is
possible to write problem \eqref{eq:darcy} in a more compact form similar to the
standard Darcy formulation. To this end, we introduce the following divergence
operator between dimensions as $\mathfrak{D} \cdot: \Omega \times \Gamma
\rightarrow \Omega$ such that
\begin{gather*}
    \mathfrak{D} \cdot \bm{w} = \nabla \cdot \check{\bm{w}}  -
    \jump{\hat{\bm{q}}\cdot\hat{\bm{n}}},
\end{gather*}
as well as a mixed-dimensional gradient operator $\mathfrak{D}: \Omega
\rightarrow \Omega \times \Gamma$ such that
\begin{gather*}
    \mathfrak{D} q = [\nabla \check{q}, \check{q} - \hat{q}]^\top.
\end{gather*}
Considering the Darcy velocity composed by $\bm{u}$ in $\Omega$ and
$\bm{u}\cdot\bm{n}$ on $\Gamma$, system \eqref{eq:darcy} becomes
\begin{gather}\label{eq:darcy_mix_form}%
    \tag{\ref{eq:darcy_d}, \ref{eq:darcy_cc}-bis}
    \begin{aligned}
        &\bm{u} + \mathfrak{K} \mathfrak{D} p = 0
        &&\quad \text{in } \Omega \times \Gamma\\
        &\mathfrak{D} \cdot \bm{u} = f&
        &\quad \text{in } \Omega
    \end{aligned}.
\end{gather}
In the previous equation $\mathfrak{K}$ is defined accordingly to include both
the effective tangential and normal permeabilities.


\subsection{Mixed-dimensional transport problem} \label{subsec:mix_transport}

Once the problem \eqref{eq:darcy} is solved, the velocity field $\bm{u}$ can be
used to describe the mixed-dimensional transport problem. We consider the same
splitting of $\Omega$ as in Subsection \ref{subsec:mix_darcy}. We introduce the
concentration compound, following the idea presented in \cite{Fumagalli2012a},
as $c=(0, c^1, \ldots, c^N, 0)$, with $0 \leq c \leq 1$ for all
$d$, to describe the portion of a scalar passive in a grid cell. We indicate by
$t$ the time variable and $(0, T)$ the time interval.

The model is the following, given $\bm{u}$ find $c$ such that
\begin{subequations} \label{eq:transport}
\begin{gather}\label{eq:transport_c}
    \check{\phi} \check{\epsilon} \partial_t \check{c} +
    \nabla \cdot ( \check{\bm{u}} \check{c}) -
    \jump{\hat{\bm{u}} \cdot \hat{\bm{n}} \hat{c}} = \check{r}
    \quad \text{in } \Omega \times (0, T),
\end{gather}
where $\epsilon$ represents the ``fracture aperture'' in each dimension, we have
$\epsilon = (\epsilon^0, \ldots, \epsilon^{N-1}, 1)$.  For fracture
intersections, the aperture is naturally interpreted as the area of the
cross-section (or volume in 0d). We assume that $\epsilon^d>0$ for all $0\leq
d\leq N$. For simplicity we assume that $\epsilon$ does
not depend on time. We have indicated with $\partial_t {c}$ the time derivative
of $c$. Finally, $\phi$ and $r$ are the porosity and a scalar source/sink term
in each dimension, represented as compounds.

In \eqref{eq:transport_c} the divergence term models
the conservation of $c$ in the current dimension $d$ and the flow exchange with
the lower dimension $d-1$. The jump operator describes the flow exchange with
the higher dimension $d+1$. The coupling conditions, related to the definition
of inflow/outflow of the lower-dimensional objects, associated to
\eqref{eq:transport_c} are
\begin{gather}\label{eq:transport_cc}
    \hat{c}|_\pm(\bm{x}) =
    \begin{dcases*}
        \hat{c}|_\pm(\bm{x}) & if $\hat{\bm{u}} \cdot
        \hat{\bm{n}}|_\pm \geq 0$\\
        \check{c}(\bm{x}) & otherwise
    \end{dcases*}
    \quad \text{on } \Gamma \times (0, T).
\end{gather}
Also in this case, problem \eqref{eq:transport} couples the concentration
defined on $\Omega^{d}$ with $\Omega^{d-1}$ and $\Omega^{d+1}$, but not directly
with $\Omega^{d-2}$ or $\Omega^{d+2}$.

With \eqref{eq:transport_c} we need to assign initial condition
$\overline{\overline{c}}$ for the concentration $c$, as
\begin{gather}\label{eq:transport_ic}
    c(t=t_{ini}) = \overline{\overline{c}} \quad \text{on } \Omega \times \{0\},
\end{gather}
with $t_{ini}$ the initial time. Finally, on the inflow part of the outher
boundary we assign a boundary condition $\overline{c}$ for \eqref{eq:transport_c}
\begin{gather} \label{eq:transport_bc}
    c = \overline{c} \quad \text{on } \partial \Omega^d_{out,
    \bm{u}\cdot\bm{n} > 0} \times (0, T),
\end{gather}
where $\partial \Omega^d_{out, \bm{u}\cdot\bm{n} > 0}$ is the portion
of $\partial \Omega^d_{out}$ such that an inflow occurs.
\end{subequations}
Equation \eqref{eq:transport} is the mixed-dimensional transport problem.

It is possible to write \eqref{eq:transport} with the formalism introduce at the
end of the previous part, obtaining
\begin{gather}\label{eq:transport_mix_form}%
    \tag{\ref{eq:transport_c}-bis}
    \phi \epsilon \partial_t c + \mathfrak{D} \cdot (\bm{u}c) = r
    \quad \text{in } \Omega \times (0, T),
\end{gather}
which is a general form of a conservative equation in the mixed-dimensional
setting.




\section{Weak and integral formulation} \label{sec:weak}


In the view of introducing the numerical approximations, in this part we
consider the weak formulation of \eqref{eq:darcy}
and the integral formulation of \eqref{eq:transport}.


\subsection{Weak formulation of mixed-dimensional Darcy problem}


We introduce the functional spaces suitable to approximate the pressure and the
velocity. Considering a
fixed-dimensional domain $\Omega^d$ we have for the
pressure $p^d \in Q^d = L^2(\Omega^d)$, for $d>0$, and $p^0 \in
Q^0 = \mathbb{R}$. For the velocity we consider, for $d>0$, the
Hilbert space
\begin{gather*}
    V^d = \left\{ \bm{v} \in H_{div}(\Omega^d): \bm{v}\cdot\bm{n}|_\pm \in
    L^2(\Gamma^{d-1})\right\}
\end{gather*}
The additional condition in $V^d$ is related to the Robin-type nature of the condition
\eqref{eq:darcy_cc}, see \cite{Martin2005,Fumagalli2016a}. The global spaces
$Q$ and $V$ are defined as the union of the local spaces.

We introduce now the bilinear forms associated to the problem \eqref{eq:darcy}.
We have $a^d(\cdot, \cdot): V^d \times V^d \rightarrow \mathbb{R}$ and
$b^d(\cdot, \cdot):
V^d \times Q^d \rightarrow \mathbb{R}$, defined as
\begin{gather*}
    a^d(\bm{w}, \bm{v}) = (K^{-1}\bm{w}, \bm{v})_{\Omega^d}
    \qquad
    b^d(\bm{w}, q) = (\nabla \cdot \bm{w}, q)_{\Omega^d},
\end{gather*}
with $\bm{w}, \bm{v} \in V^d$ and $q \in Q^d$.
For the coupling between dimensions, $d>0$, we introduce
$\alpha^{d}(\cdot, \cdot): V^{d} \times V^{d} \rightarrow \mathbb{R}$ and
$\beta^{d, d-1}(\cdot, \cdot): V^{d} \times Q^{d-1} \rightarrow \mathbb{R}$
defined as
\begin{gather*}
    \alpha^{d}(\bm{w}, \bm{v}) = \sum_{\pm}(\eta_*^d \bm{w} \cdot \bm{n}|_{\pm},
    \bm{v} \cdot \bm{n}|_{\pm})_{\Omega^{d-1}}\\
    \beta^{d, d-1}(\bm{w}, q) = \sum_{\pm} (\bm{w} \cdot \bm{n}
    |_{\pm}, q)_{\Omega^{d-1}},
\end{gather*}
where $\bm{w}, \bm{v} \in V^d$ and $q \in Q^{d-1}$. In the previous definition
we have $\eta_*^d$ the inverse of $k_*^d$ and $\bm{n}$ the normal of
$\Omega^{d-1}$ in the tangent space of $\Omega^d$.
The global bilinear forms are defined by $a(\cdot, \cdot): V \times V \rightarrow \mathbb{R}$
and by $b(\cdot, \cdot): V \times Q \rightarrow \mathbb{R}$ as
\begin{gather*}
    a(\bm{w}, \bm{v}) = \sum_{d>0} a^d(\bm{w}^d, \bm{v}^d) +
    \alpha^{d}(\bm{w}^d, \bm{v}^d)\\
    b(\bm{w}, q) = \sum_{d>0} b^d(\bm{w}^d, q^d) + \beta^{d, d-1}(\bm{w}^d,
    q^{d-1}),
\end{gather*}
with $\bm{w}, \bm{v} \in V$ and $q \in Q$.
We can introduce the weak formulation of \eqref{eq:darcy}, which reads:
find $(\bm{u}, p) \in V \times Q$ such that
\begin{align}\label{eq:darcy_weak}
    \begin{aligned}
        &a(\bm{u}, \bm{v}) + b(\bm{v}, p) = G(\bm{v}) & \forall \bm{v} \in V\\
        &b(\bm{u}, q) =  F(q) & \forall q \in Q
    \end{aligned}
\end{align}
where the functionals are defined as
\begin{gather*}
    F(q) = \sum_d (f^d, q^d)_{\Omega^d} \quad
    G(\bm{v}) = -\sum_{d>0} (\bm{v} \cdot \bm{n} |_{\partial
    \Omega^d_{out}}, \overline{p}^d)_{\partial \Omega^d_{out}}.
\end{gather*}


\subsection{Integral formulation of mixed-dimensional transport problem}


In Subsection \ref{subsec:fv} we consider a finite volume approximation of
problem \eqref{eq:transport}, thus we introduce here its integral formulation.

Considering a suitable sub-domain of $\Omega^d$ called $E$, later will be a cell of the
computational grid, we integrate equation \eqref{eq:transport} obtaining
the following
\begin{gather} \label{eq:transport_int}
    \begin{gathered}
        (\check{\phi}\check{\epsilon}\partial_t \check{c}, 1)_E + (
        \check{\bm{u}} \cdot \check{\bm{n}}|_{\partial E}, \check{c}
        |_{\partial E})_{\partial E} +\\
        \sum_\pm (\check{\bm{u}} \cdot
        \check{\bm{n}}|_{I, \pm}, \check{c}|_{I, \pm})_{I, \pm} -
        ( \hat{\bm{u}} \cdot \hat{\bm{n}}|_{E, \pm}, \hat{c}|_{E, \pm})_{E} =
        (\check{r}, 1)_E
    \end{gathered}
\end{gather}
where $I = E \cap \Omega^{d-1}$ is the intersection between $E$ and the lower
dimensional domain $\Omega^{d-1}$, which can be an empty set.

In \eqref{eq:transport_int}, we notice the term depending on $I$ which describe
the inflow/outflow exchange between the current and lower dimensions.
Condition \eqref{eq:transport_cc} applies also in this case.





\section{Numerical discretization} \label{sec:numerical_discretization}


In this section we introduce the numerical scheme used to solve the
mixed-dimensional problems \eqref{eq:darcy} and \eqref{eq:transport}. In
particular, for the former we consider the lowest order dual virtual element
method approximation (VEM), while for the latter we employ a finite volume
discretization with an upwind scheme and implicit Euler in time.

Let us introduce a suitable tessellation $\mathcal{T}$ of $\Omega$ in
polytopes indicated by $E$, such that
\begin{gather*}
    \mathcal{T} = \cup_i E_i
    \quad \text{and} \quad
    E_i \cap E_j = \emptyset \text{ for } i \neq j.
\end{gather*}
We denote by $h_E$ the diameter of $E$, by $\bm{x}_E$ the centre of $E$, by
$\mathcal{E}(E) = \left\{ e \in \partial E \right\}$ the set of faces for the
cell $E$, and by $\bm{n}_e$ with $e
\in \mathcal{E}(E)$ the unit normal of $e$ pointing outward with respect to the
internal part of $E$. The term faces indicates proper faces for the 3d grid,
edges for the 2d grids, and vertices for the 1d grids.
We indicate also by $h = \max_{E \in \mathcal{T}} h_E$ the
characteristic grid size and by $\mathcal{E} = \cup_{E \in \mathcal{T}}
\mathcal{E}(E)$ the set of the grid faces. It is important to note that the polytopes
may have different spatial dimension, \textit{i.e.} they are 3d or 2d or 1d or 0d
objects. However, fixing a single dimension $d$ the cells of the grid belongs
to $\Omega^d$.

For simplicity, both problems \eqref{eq:darcy} and \eqref{eq:transport} use the
same computational grid.

In Subsection \ref{subsec:dual_vem} we present the discretization for problem \eqref{eq:darcy},
while in Subsection \ref{subsec:fv} the numerical approximation for problem
\eqref{eq:transport}. Finally, in Subsection \ref{subsec:coarsening} we
introduce a coarsening strategy used in the numerical examples.


\subsection{Mixed-dimensional dual VEM discretization}%
\label{subsec:dual_vem}


We present a numerical scheme to solve the mixed-dimensional Darcy problem
\eqref{eq:darcy} in presence of
polytopes in the grid. The resulting scheme will be locally and globally conservative, thus
suitable to approximate the velocity field used in \eqref{eq:transport}. We
consider the dual virtual element method of lowest order degree, for more
details on the derivation and on the analysis see
\cite{BeiraodaVeiga2014a,Brezzi2014,BeiraodaVeiga2014b,Benedetto2014,BeiraoVeiga2016,Benedetto2016,Fumagalli2016a}.

\begin{figure}[htb]
    \centering
    \includegraphics[width=0.375\textwidth]{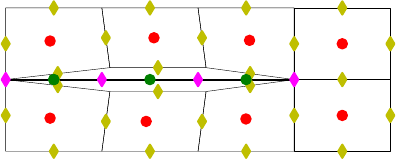}%
    \caption{Representation of the degrees of freedom for a 2d and 1d grid. The
    pressure \textit{dof} are represented by circles, red for the 2d grid and
    green for the 1d. The velocity \textit{dof} are depicted by yellow diamonds
    for the 2d grid and purple diamonds for the 1d. The nodes of the 2d grid are moved
    only for visualization purpose.}%
    \label{fig:split_dot}%
\end{figure}

We introduce finite dimensional spaces to approximate the Darcy velocity
$\bm{u}$ and the pressure $p$ in each element of the grid. For simplicity we
consider a single dimension and, if not essential, we drop the superscript $d$ to simplify
notation. Again the differential operators are understood to be defined on the
tangential space. Given an element $E$ the local discretization space for the
pressure is $Q_h^d(E) = \mathbb{P}_0(E) \subset L^2(E)$, where $\mathbb{P}_r(E)$
is the space of polynomial of order $r$ on the domain $E$. For the velocity we
need to introduce the following space
\begin{gather*}
    V_h^d(E) = \left\{ \bm{v} \in H_{div}(E):\, \bm{v}\cdot\bm{n}_e \in
    \mathbb{P}_0(e) \forall e \in \mathcal{E}(E), \right. \\
    \nabla \cdot \bm{v} \in
    \left. \mathbb{P}_0(E), \nabla \times \bm{v} = \bm{0} \right\}.
\end{gather*}
The shape of the functions in $V_h^d(E)$ is not defined a-priori and are
implicitly defined by $V_h^d(E)$.  The curl-free conditions necessary to
uniquely define the elements in $V_h^d(E)$. We indicate by $V_h^d$ the velocity
approximation space in the same dimensional grid and by $V_h$ the global
discretization space for the velocity formed by the compound
$(V_h^1,...,V_h^N)$. Similarly, $Q_h^d$ indicates the pressure approximation
space in the same dimension $d$ while $Q_h$ is the global discretization space.
For the velocity, we impose to the faces which are not in contact with different
dimensions to be single value.  Otherwise, the degree of freedom is doubled and
connected thorough the coupling condition \eqref{eq:darcy_cc} to the lower
dimensional object. See Figure \ref{fig:split_dot} for an example.

With the previous definition of $Q_h$ and $V_h$ it is immediate to approximate
the bilinear form $b$ for all the dimensions as well as the bilinear forms
associated with the coupling conditions, $\alpha^d$ and $\beta^{d, d-1}$. The
functionals $F$ and $G$ are discretized similarly. The bilinear forms $a^d$ are
not immediately computable with the degrees of freedom introduced, but require
the definition of a suitable projection operator.
First we introduce the local space
\begin{gather*}
    \mathcal{V}_h^d(E) = \left\{ \bm{v} \in V_h^d(E): \bm{v} = K \nabla v, \text{for a
    } v \in \mathbb{P}_1(E)\right\},
\end{gather*}
and the projection operator is defined as $\Pi_0: V_h^d(E) \rightarrow
\mathcal{V}_h^d(E)$ such that given $\bm{v} \in V_h^d(E)$ we have
$a^d(\bm{v} - \Pi_0 \bm{v}, \bm{w}) = 0$ for all $\bm{w} \in \mathcal{V}_h^d(E)$.
Note that the space $\mathbb{P}_1(E)$ in the definition of $\mathcal{V}_h^d(E)$
will be approximated by a (tangential) monomial basis.
With this property it is possible to split the bilinear form $a^d$ in a term
on $\mathcal{V}_h^d(E)$ and one on the $a^d$-orthogonal space of
$\mathcal{V}_h^d(E)$, namely
\begin{gather*}
    a^d(\bm{u}, \bm{v}) = a^d(\Pi_0 \bm{u}, \Pi_0 \bm{v}) +
    a^d(T_0 \bm{u}, T_0 \bm{v})
\end{gather*}
with $T_0 = I - \Pi_0$. The first bilinear form is now fully computable with the
velocity degrees of freedom introduced before and represent a consistency term
with respect to $a^d(\bm{u}, \bm{v})$. The second term is not yet computable and
can be replaced by a stabilization term. Following the ideas presented in
\cite{Brezzi2014,BeiraodaVeiga2014b,BeiraoVeiga2016}, we approximate by
\begin{gather*}
    a^d(T_0 \bm{u}, T_0 \bm{v}) \approx \varsigma s^d(\bm{u}, \bm{v})
\end{gather*}
where $s^d: V_h^d(E)\times V_h^d(E)\rightarrow \mathbb{R}$ is the bilinear form associated
with the stabilization and $\varsigma = \varsigma(d) \in \mathbb{R}^+$ is a scaling parameter
described in the sequel. In details, denoting by $\bm{\varphi}$ an element of
the basis for $V_h^d(E)$, in our case we consider
\begin{gather*}
    s^d(\bm{\varphi}_\omega, \bm{\varphi}_\theta) = \sum_{i=1}^{N_{dof}} dof_i(T_0
    \bm{\varphi}_\omega) \, dof_i(T_0 \bm{\varphi}_\theta),
\end{gather*}
with $N_{dof}$ the total number of velocity degrees of freedom of the element $E$.
We introduce the discrete version $a^d_h: V_h^d(E) \times V_h^d(E) \rightarrow
\mathbb{R}$ of the bilinear form $a^d$, defined as
\begin{gather*}
    a_h^d(\bm{u}, \bm{v}) = a^d(\Pi_0 \bm{u}, \Pi_0 \bm{v}) +
    \varsigma s^d(\bm{u}, \bm{v}),
\end{gather*}
and the discrete version of the weak problem \eqref{eq:darcy_weak}, which reads:
find $(\bm{u}, p) \in V_h\times Q_h$ such that
\begin{align}\label{eq:darcy_weak_discr}
    \begin{aligned}
        &a_h(\bm{u}, \bm{v}) + b(\bm{v}, p) = G(\bm{v}) & \forall \bm{v} \in V_h\\
        &b(\bm{u}, q) =  F(q) & \forall q \in Q_h
    \end{aligned}
\end{align}
where $a_h(\bm{u}, \bm{v}) = \sum_{d>0} a_h^d(\bm{u}^d, \bm{v}^d)
+\alpha^{d}(\bm{w}^d, \bm{v}^d)$.
How to represent the discrete problem in term of local matrix computations
refer to \cite{BeiraodaVeiga2014a,Fumagalli2016a}.

The stabilization parameter $\varsigma$ is used to impose the scaling on $h$
of the stabilization
term equivalent to the consistency term. In practice we require that,
for a fixed dimension $d$, there exists $\iota_*, \iota^* \in \mathbb{R}^+$,
independent from the discretization, such that
\begin{gather*}
    \iota_* a^d(\Pi_0 \bm{v}, \Pi_0 \bm{v}) \leq \varsigma s^d(\bm{v}, \bm{v}) \leq
    \iota^* a^d(\Pi_0 \bm{v}, \Pi_0 \bm{v}) \quad \forall \, \bm{v} \in V_h^d.
\end{gather*}
Following \cite{Fumagalli2016a}, by a scaling computation it is possible to
evaluate the dependency of $\varsigma$ from the local grid size $h_E$, obtaining
the following relation for a fix dimension $d$
\begin{gather*}
    \varsigma(E) = h_E^{2-d}.
\end{gather*}
Note that this expression is local and independent from the maximal dimension
$N$ of the problem.

We finally introduce the matrix formulation associated with the problem
\eqref{eq:darcy_weak_discr}. Considering the following matrices
\begin{gather*}
    [A_d]_{ij} = a^d_h(\bm{\varphi}_j, \bm{\varphi}_i) +
    \alpha^d(\bm{\varphi}_j, \bm{\varphi}_i),
    \\
    [B_d]_{ij} = b^d(\bm{\varphi}_j, \psi_i),
    \quad
    [C_{d, d-1}]_{ij} = \beta^{d, d-1}(\bm{\varphi}_j, \psi_i)
\end{gather*}
and vectors
\begin{gather*}
    [G_d]_i = G(\bm{\varphi}_i), \quad
    [F_d]_i = F(\psi_i),
\end{gather*}
where $\phi$ are basis for the pressure and $[\cdot]_{i,j}$ indicates the
element $(i, j)$ in the matrix, a similar notation is used for vectors.
The global problem reads for $N=3$, solve the following linear system
\begin{gather*}
    \label{eq:linsys}
    \begin{bmatrix}
        A_3             & B_3 & 0        & C_{3, 2} & 0 & 0 & 0 & 0\\
        B_3^\top        &  0   &  0       &  0 & 0 & 0 & 0 & 0\\
        0 & 0 & A_2             & B_2 & 0        & C_{2, 1} & 0 & 0\\
        C_{3, 2}^\top & 0 & B_2^\top        &  0   &  0       &  0 & 0 & 0\\
        0 & 0 & 0               &  0   & A_{1}      & B_{1} & 0 & C_{1, 0}\\
        0 & 0 & C_{2, 1}^\top &   0  & B_{1}^\top & 0 & 0 & 0 \\
        0 & 0 & 0 & 0 & 0 & 0& I &  0 \\
        0 & 0 & 0 & 0 & C_{1, 0}^\top & 0 & 0 & 0
    \end{bmatrix}
    \begin{bmatrix}
        \bm{u}_3\\
        p_3\\
        \bm{u}_2\\
        p_2\\
        \bm{u}_1\\
        p_1\\
        \bm{u}_0\\
        p_0
    \end{bmatrix} =
    \begin{bmatrix}
        G_3\\
        F_3\\
        G_2\\
        F_2\\
        G_1\\
        F_1\\
        0\\
        0
    \end{bmatrix}
\end{gather*}
Where $\bm{u}_d$ and $p_d$ represent the vectors associated with the degrees of
freedom of velocity and pressure in each dimension $d$.  We notice that the
$\bm{u}_0$ is considered only for clearness and to preserve the structure of the
matrix. In practice, it is possible to remove $\bm{u}_0$ from the system.


\subsection{Mixed-dimensional FV discretization}%
\label{subsec:fv}

\begin{figure*}[!t]
    \centering
    \resizebox{0.3\textwidth}{!}{%
    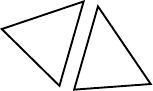}%
    \resizebox{0.3\textwidth}{!}{%
    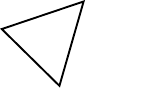}%
    \resizebox{0.34\textwidth}{!}{%
    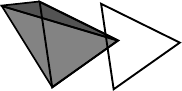}%
    \caption{Representation of the coupling between dimension for the upwind
        discretization scheme. On the left between two cells in the same
        dimension $d$, in the centre between a $d$-dimensional cell and a
        $d-1$-dimensional cell, and on
        the right between a $d+1$-dimensional cell and a $d$-dimensional cell.}%
    \label{fig:upwind}%

    \vspace*{\floatsep}

    \centering
    \includegraphics[width=0.3\textwidth]{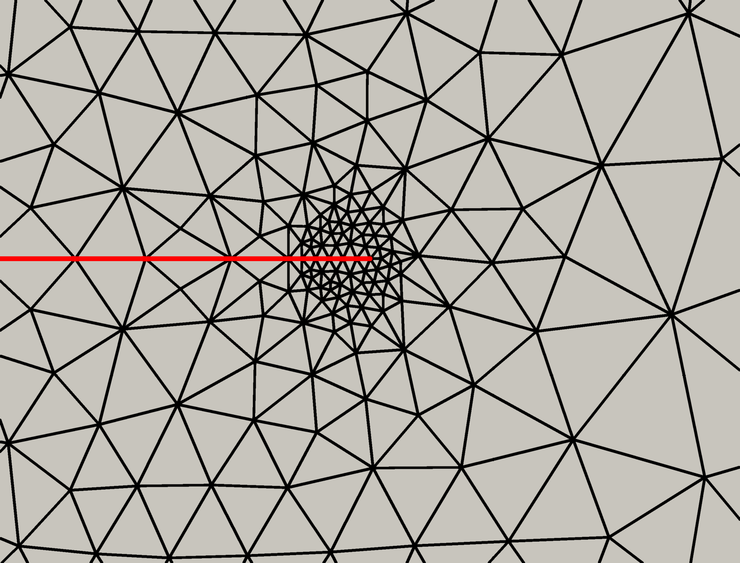}\hspace{0.15\textwidth}%
    \includegraphics[width=0.3\textwidth]{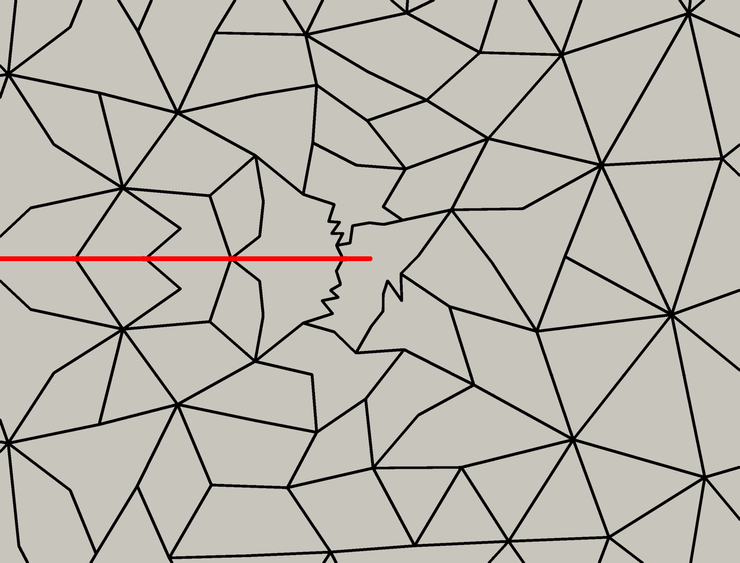}
    \caption{Example of the coarsening strategy adopted. In the left part the
    computational grid is artificially forced to be more fine at the tip of a
    fracture, while in the right the resulting grid after the coarsening.
    Clearly the cell measures are comparable.}%
    \label{fig:coarsening_algo}
\end{figure*}

We consider now the discretization of equation \eqref{eq:transport_int}. For
simplicity, we consider a fixed time step $\Delta t$ such that $T/\Delta t$ is
an integer number. An implicit Euler is consider in time obtaining the
semi-discrete version of \eqref{eq:transport_int}
\begin{gather*}
    \begin{gathered}
        \dfrac{(\check{\phi}\check{\epsilon}(\check{c}^{n+1} - \check{c}^{n}),
        1)_E}{\Delta t} + (
        \check{\bm{u}} \cdot \check{\bm{n}}|_{\partial E}, \check{c}^{n+1}
        |_{\partial E})_{\partial E} +\\
        \sum_\pm (\check{\bm{u}} \cdot
        \check{\bm{n}}|_{I, \pm}, \check{c}^{n+1}|_{I, \pm})_{I, \pm} -
        ( \hat{\bm{u}} \cdot \hat{\bm{n}}|_{E, \pm}, \hat{c}^{n+1}|_{E, \pm})_{E} =
        (\check{r}^{n+1}, 1)_E,
    \end{gathered}
\end{gather*}
where the superscript indicates the time step index, clearly $c^0 =
\overline{\overline{c}}$. The approximation of the boundary integrals in the
previous equation relay on an upwind scheme. We introduce a Kronecker-type delta
as
\begin{gather*}
    \delta_{\bm{u}\cdot\bm{n}} =
    \begin{dcases*}
        1 & if $\bm{u}\cdot\bm{n} \geq 0$\\
        0 & else
    \end{dcases*}.
\end{gather*}
Given a cell face $e \in \mathcal{E}(E)$ we have
\begin{gather*}
    (\check{\bm{u}} \cdot \check{\bm{n}}|_{e},
    \check{c}^{n+1}|_{e})_{e} = \check{\bm{u}} \cdot \check{\bm{n}}|_{e} [
    \delta_{\check{\bm{u}} \cdot \check{\bm{n}}|_{e}} \check{c}^{n+1}(E) +
    (1-{\delta}_{\check{\bm{u}} \cdot \check{\bm{n}}|_{e}}) \check{c}^{n+1}(L)]
\end{gather*}
where the cells $K$ and $L$ share the face $e$.
Given a cell $E$ such that one of its face $e$ intersect one side of $I$, we get
\begin{gather*}
    (\hat{\bm{u}} \cdot \hat{\bm{n}}|_{e},
    \hat{c}^{n+1}|_{e})_{e} =
    \hat{\bm{u}} \cdot \hat{\bm{n}}|_{e} [ \delta_{\hat{\bm{u}} \cdot
    \hat{\bm{n}}|_{e}} \hat{c}^{n+1}(E) + (1-{\delta}_{\hat{\bm{u}}
    \cdot\hat{\bm{n}}|_{e}}) \check{c}^{n+1}(\check{e})]
\end{gather*}
where $\check{e}$ indicates a cell in the co-dimensional grid which is in
communication with the face $e$.
Finally, on side of the last term of the semi-discrete problem can be approximated by
\begin{gather*}
    (\hat{\bm{u}} \cdot \hat{\bm{n}}|_{E},
    \hat{c}^{n+1}|_{E})_{E} =
    \hat{\bm{u}} \cdot \hat{\bm{n}}|_{E}[\delta_{\hat{\bm{u}} \cdot
    \hat{\bm{n}}|_{e}} \hat{c}^{n+1}(\hat{E}) + (1-
    {\delta}_{\hat{\bm{u}}\cdot\hat{\bm{n}}|_{e}}) \check{c}^{n+1}(E)]
\end{gather*}
where $\hat{E}$ represent the cell in the higher dimensional grid which has a
face in communication with the cell $E$.
The previous condition applies to both the equi and co-dimensional coupling, see
Figure \ref{fig:upwind} as an example. Note that the chosen discretization is compatible
with the coupling condition \eqref{eq:transport_cc}.

The matrix formulation of the transport problem consider the introduction of the
following matrices
\begin{gather*}
    [M_d]_{ii} = \frac{(\phi^d \epsilon^d, 1)_{i}}{\Delta t},
    \quad
    [U_d]_{ii} = \sum_{j \in \mathcal{N}(i), e \in \overline{i} \cap
    \overline{j}}
    (\bm{u}^d \cdot \bm{n}^d |_{e}, \delta_{\bm{u}^d \cdot
    \bm{n}^d|_{e}})_{e},\\
    [U_d]_{ij} = \sum_{j \in \mathcal{N}(i), e \in \overline{i} \cap
    \overline{j}}
    (\bm{u}^d \cdot \bm{n}^d |_{e},
    1-{\delta}_{\bm{u}^d \cdot\bm{n}^d|_{e}})_{e}
\end{gather*}
where $\mathcal{N}(E)$ is the set of all neighbour cells of $E$,
$i$ and $j$ indicate the generic cell $i$ and $j$, and $e$ is a face in common
with $i$ and $j$. For the coupling between dimensions we have
\begin{gather*}
    [U_{d,d-1}]_{ii} = \sum_{j \in \mathcal{N}^{d-1}(i)}
    (\bm{u}^d \cdot \bm{n}^d|_j, \delta_{\bm{u}^d \cdot \bm{n}^d|_j})_j\\
    [U_{d,d-1}]_{ij} = \sum_{j \in \mathcal{N}^{d-1}(i)}
    (\bm{u}^d \cdot \bm{n}^d|_j, 1-\delta_{\bm{u}^d \cdot \bm{n}^d|_j})_j,
\end{gather*}
where $\mathcal{N}^{d-1}(E)$ is the set of all neighbour cells of $E$ in the
lower dimensional grid. We obtain also
\begin{gather*}
    [U_{d,d+1}]_{ii} = -\sum_{j \in \mathcal{N}^{d+1}(i)}
    (\bm{u}^{d+1} \cdot \bm{n}^{d+1}|_i,
    \delta_{\bm{u}^{d+1} \cdot \bm{n}^{d+1}|_i})_i\\
    [U_{d,d+1}]_{ij} = -\sum_{j \in \mathcal{N}^{d+1}(i)}
    (\bm{u}^{d+1} \cdot \bm{n}^{d+1}|_i, 1-\delta_{\bm{u}^{d+1} \cdot
    \bm{n}^{d+1}|_i})_i,
\end{gather*}
where $\mathcal{N}^{d+1}(E)$ is the set of all neighbour cells of $E$ in the
higher dimensional grid.
We finally obtain the following linear system to be inverted
\begin{gather*}
    \begin{bmatrix}
        U_3 + M_3     & U_{3,2} & 0       & 0       \\
        U_{2,3} & U_2 + M_2     & U_{2,1} & 0       \\
        0       & U_{1,2} & U_1 + M_1     & U_{1,0} \\
        0       & 0       & U_{0,1} & 0
    \end{bmatrix}
    \begin{bmatrix}
        c_3^{n+1}\\
        c_2^{n+1}\\
        c_1^{n+1}\\
        c_0^{n+1}
    \end{bmatrix} =
    \begin{bmatrix}
        \zeta_3^{n,n+1}\\
        \zeta_2^{n,n+1}\\
        \zeta_1^{n,n+1}\\
        0
    \end{bmatrix}
\end{gather*}
Where $c_d$ represents the vector associated with the degrees of freedoms of the
concentration in each dimension $d$.
In the previous linear system $\zeta_d$
represent the right-hand side, combination of the source term $r_d^{n+1}$ and
the concentration at the previous time step. We get
\begin{gather*}
    \zeta_d^{n,n+1}= M_d c_d^{n} + r_d^{n+1}.
\end{gather*}


\subsection{Coarsening strategy}\label{subsec:coarsening}


The creation of conforming grids in presence of several fractures can be a
challenging task, especially in 3d. In this work we rely on the Gmsh mesh
generator \cite{Geuzaine2009} for the construction of the computational grids.
In presence of almost intersecting fractures or small fracture branches, the
grid may result composed by a high number of simplex cells. To overcome
this difficulty we exploit one of the advantages of the VEM able to handle cells
of arbitrary geometry. To be specific, the theory developed in
\cite{BeiraodaVeiga2014a,Brezzi2014,BeiraodaVeiga2014b,BeiraoVeiga2016} requires
star-shaped cells, however the study carry out in \cite{Fumagalli2017b} shows
that the VEM are able to handle also cells with cuts.

Motivated by these observations we introduce a coarsening scheme that merges
small simplex cells into larger polygons or polyhedra.  Starting from a given
simplex grid, the algorithm computes the measure (area or volume) of the cells.
Given the cell $c$ with the smallest measure, it will be merged to one or more
neighbouring cells, based on their respective measure, creating a new coarse
cell. The clustering stops when a cell measure threshold is reached.


The algorithm does not guarantee any regularity of the final grid, see an
example in Figure \ref{fig:coarsening_algo}.




\section{Applicative examples} \label{sec:examples}


We present some examples and test cases to asses the previous models and
numerical schemes. The first example, presented in Subsection \ref{example_1},
consider an extensive validation of the
mixed-dimensional Darcy problem solved by VEM through a
benchmark study presented in \cite{Flemisch2016a}. The second test case, in
Subsection \ref{example_2}, consider
the transport problem and analyse the impact of the coarsening strategy on the
results. Finally, in Subsection \ref{example_3} a realistic 3d example is
introduced and studied. In all the forthcoming examples we assume unitary
porosity and zero source term for the concentration and Darcy equations. The
other parameters will be specified.

In all of the forthcoming examples a ``Blue to Red Rainbow'' colour map is used.

The examples are part of the PorePy package, which is a simulation tool for
fractured and deformable porous media written in Python. See
\url{github.com/pmgbergen/porepy} for more details.


\subsection{Benchmark comparison}\label{example_1}

To validate the presented model, we consider benchmark 1 and 4
presented in \cite{Flemisch2016a}. We carry out the same type of analysis in the same
setting, for more detail of problem setting refer to the aforementioned work.


\subsubsection{Benchmark 1: Regular Fracture Network}

\begin{figure*}[p]
    \centering
    \includegraphics[width=0.33\textwidth]{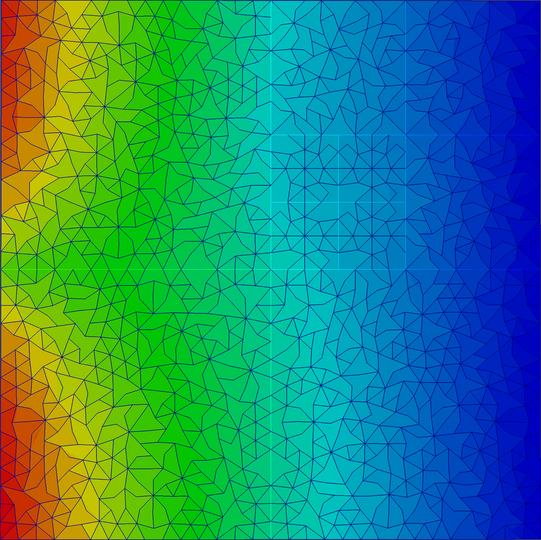}%
    \hfill%
    \includegraphics[width=0.33\textwidth]{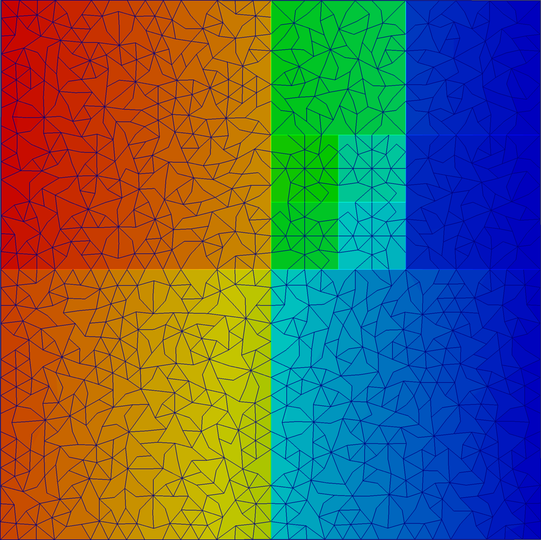}%
    \hfill%
    \includegraphics[width=0.33\textwidth]{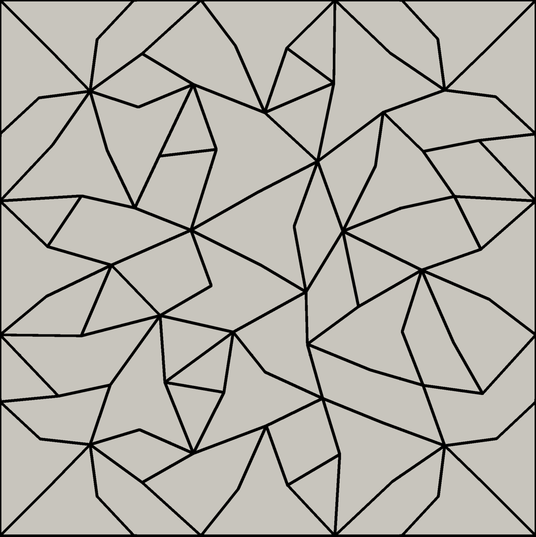}%
    \caption{Benchmark 1. On the left: pressure solution (range $(1, 1.6)$)
    with conductive fractures and computational grid. On the centre: pressure
    solution (range $(1, 3.6)$)
    with blocking fractures and computational grid. On the right: zoom of the
    grid used.}%
    \label{fig:geiger_vem}

    \vspace*{\floatsep}

    \includegraphics[width=0.33\textwidth]{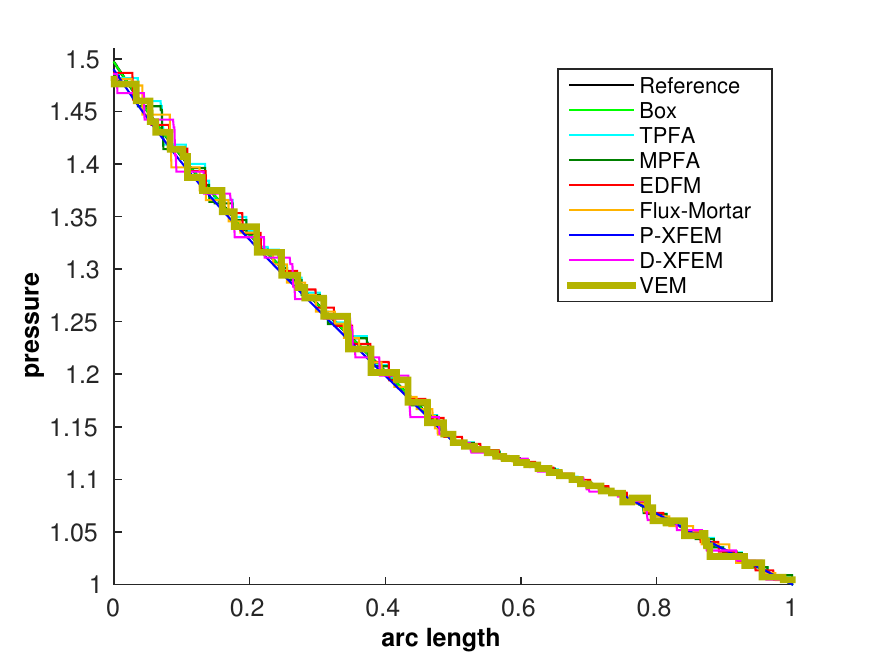}\hfill%
    \includegraphics[width=0.33\textwidth]{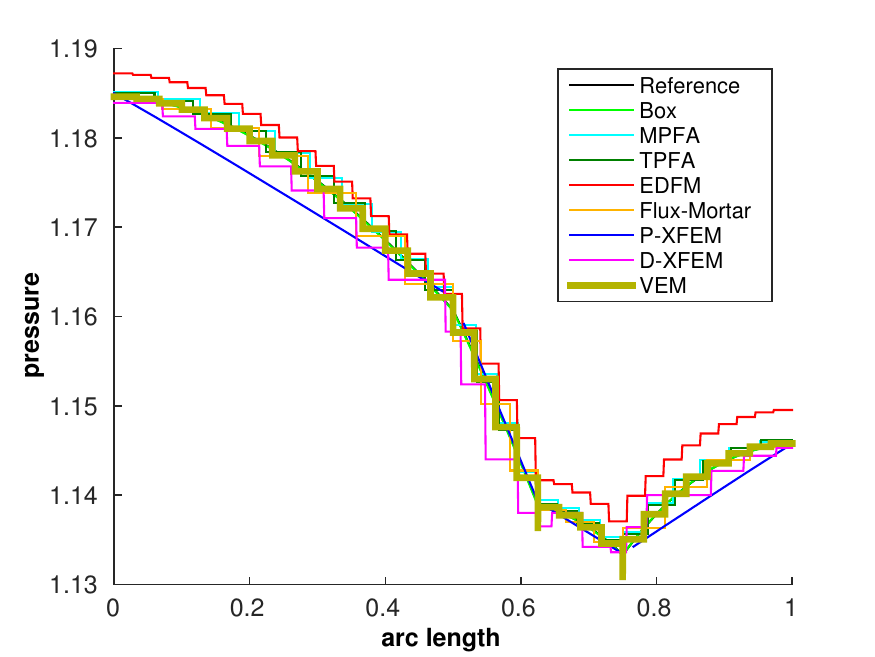}
    \includegraphics[width=0.33\textwidth]{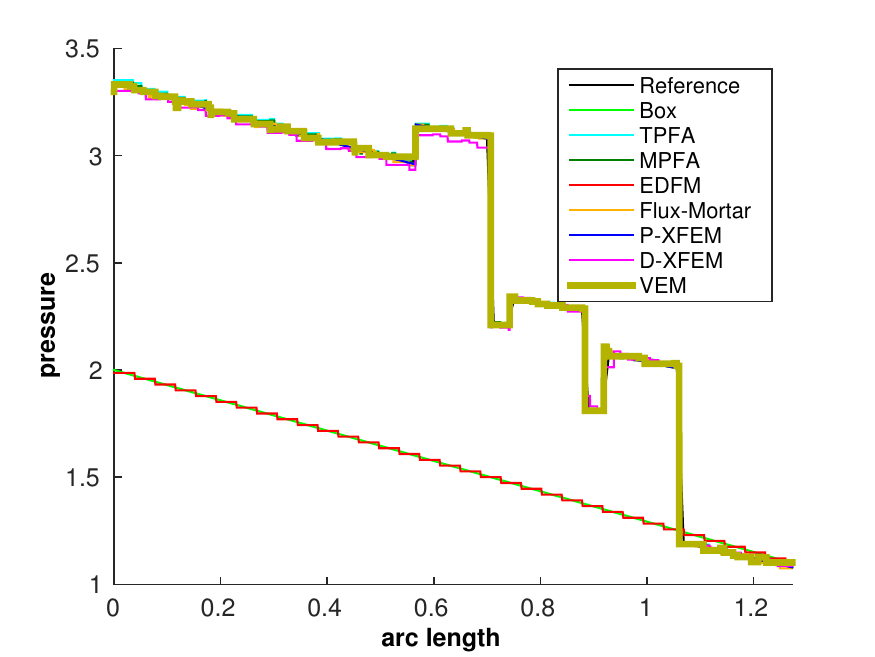}
    \caption{Benchmark 1 with conductive fractures.
    On the left: pressure along horizontal line at $y=0.7$ with permeable
    fracture. On the centre: pressure along vertical fracture at $x=0.5$ with
    permeable fractures. On the right: pressure along the line $(0.0,
    0.1)-(0.9, 1.0)$ with blocking fractures.}%
    \label{fig:geiger_over_line}

    \vspace*{\floatsep}

    \centering
    \includegraphics{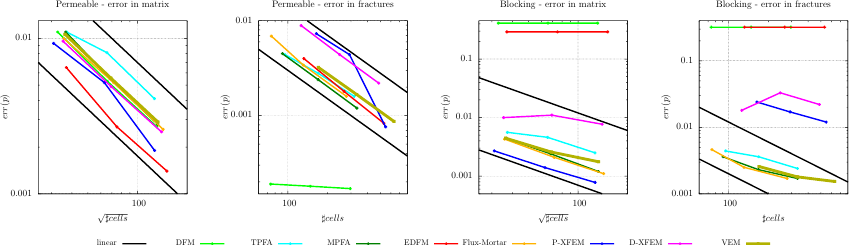}
    \caption{Benchmark 1, error evolution for the matrix and the fractures
    with conductive and blocking fractures.}%
    \label{fig:geiger_permeable_error}
\end{figure*}

The problem is inspired by \cite{Geiger2011} with different boundary conditions and
material properties. We have unitary permeability matrix and fracture aperture
equal to $10^{-4}$. We consider two possibility for fracture permeability: high
conductive with permeability $10^{4}$ and low conductive with permeability
$10^{-4}$. In the former case the solution obtained with the method presented
previously along with the computational grid are presented in Figure
\ref{fig:geiger_vem} left, the latter in the right.
We can notice a good agreement between the computed and the reference
solutions, reported in \cite{Flemisch2016a}. We point out that some of the elements present in Figure
\ref{fig:geiger_vem} are not convex.

To have a more detailed comparison we consider two plots over line for the
permeable case and one for the blocking case, shown in Figure
\ref{fig:geiger_over_line}. From now on, the method presented in this paper is
labelled as VEM.
Also in this case we have a good agreement between the computed and reference
solutions, comparable to others methods able to represent blocking fractures.
Small oscillations are related to mesh effects.

Finally, in both cases we consider the error decay for both the rock matrix
and the system of fractures. We consider a family of three meshes where the
coarsening is applied in all the cases. We obtain again non-convex elements in all
the meshes. Figure \ref{fig:geiger_permeable_error}
plots for conductive and blocking fractures the error decay.
The errors are comparable with others methods able to represent permeable and
blocking fractures. In the latter case we notice a stagnancy of the fracture
error bounding the order of converge, this phenomena is common also for other
methods presented in the benchmark study \cite{Flemisch2016a}.


\subsubsection{Benchmark 4: a Realistic Case}\label{subsub:sotra}

We consider now a complex system of 64 fractures from a real outcrop. We
consider constant rock permeability equal to $10^{-14}$m$^2$, uniform fracture
permeability $10^{-8}$m$^2$, and fracture aperture $10^{-2}$m. We impose a
pressure gradient at the boundary from the left (1013250 Pa) to the right (0
Pa). Also in this case, a triangular grid is coarsened to reduce the grid complexity.
The reference grid is composed by 12472 2d cells, 1317 1d cells, and 85 0d
cells. The coarse algorithm decreases the 2d cells down to 4703. The total
number of degrees of freedom is 19075 for the coarse grid.

\begin{figure*}[p]
    \centering
    \includegraphics[width=0.33\textwidth]{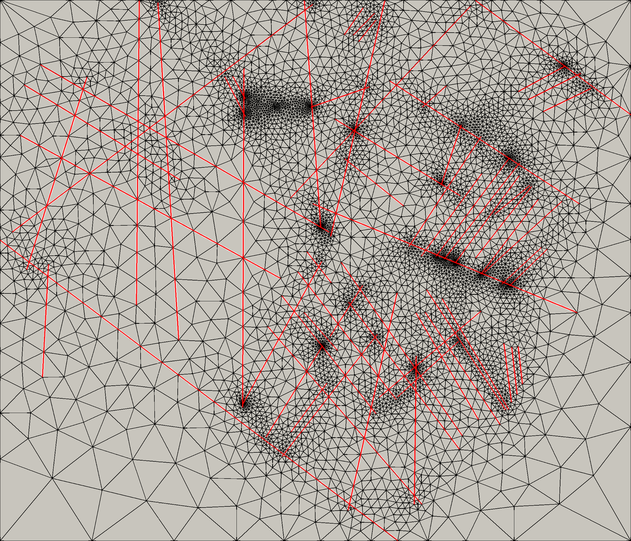}%
    \hspace*{0.02\textwidth}%
    \includegraphics[width=0.33\textwidth]{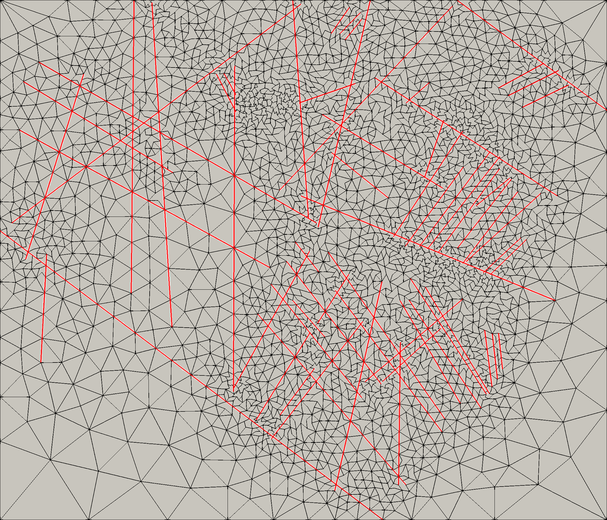}\\
    \includegraphics[width=0.33\textwidth]{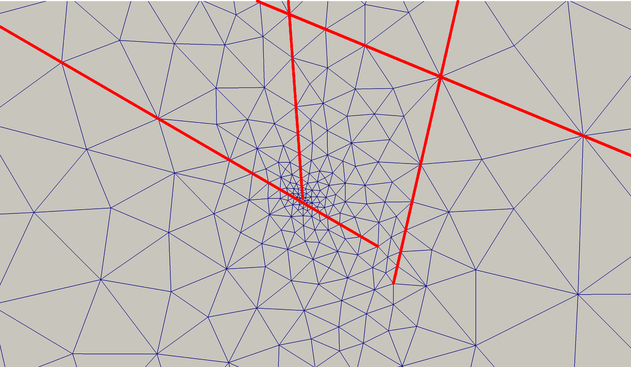}%
    \hspace*{0.02\textwidth}%
    \includegraphics[width=0.33\textwidth]{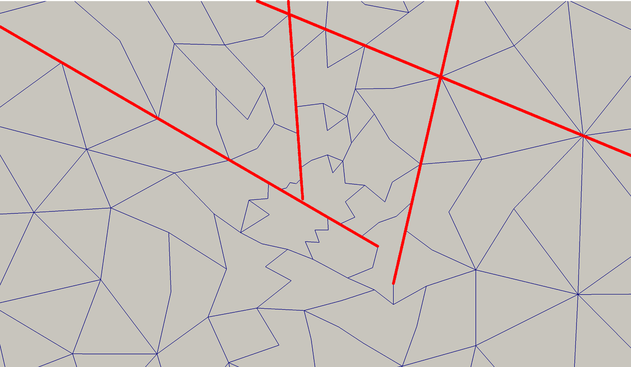}%
    \caption{Benchmark 4. On the top: original grid composed by 12302 2d-cells
    and 35153 VEM \textit{dof} and coarsened grid composed by 4599 2d-cells
    and 18803 VEM \textit{dof}. On the bottom: a zoom on almost intersecting
    fractures for the original and coarsened grids, respectively. The zoom
    is referred to the small rectangle at position,
    approximately, $(360, 350)$.}%
    \label{fig:sotra_grid}

    \vspace*{\floatsep}

    \centering
    \includegraphics[width=0.33\textwidth]{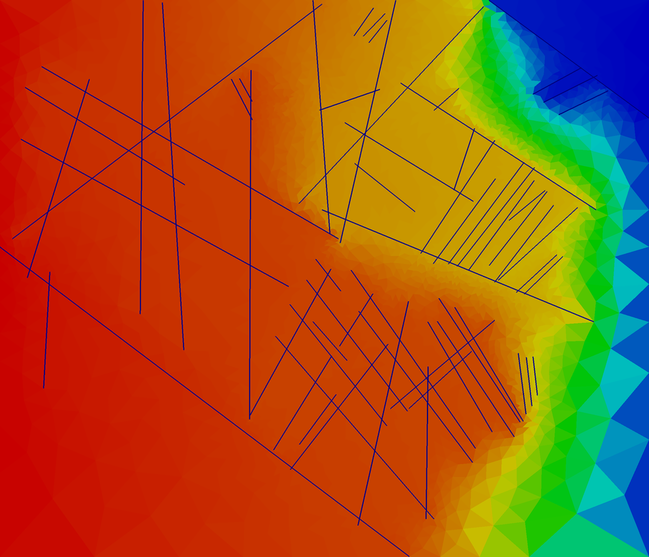}%
    \hfill%
    \includegraphics[width=0.33\textwidth]{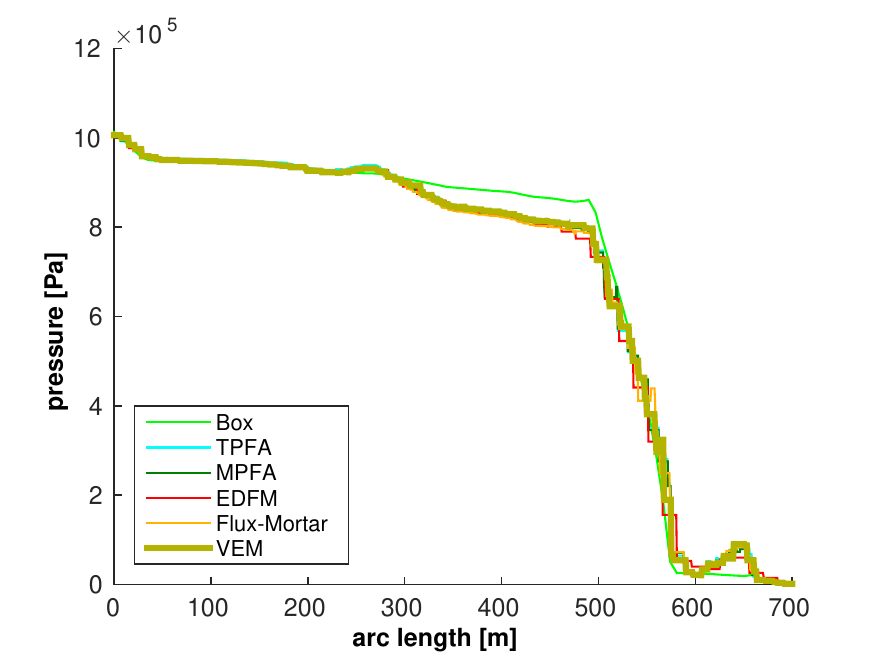}%
    \hfill%
    \includegraphics[width=0.33\textwidth]{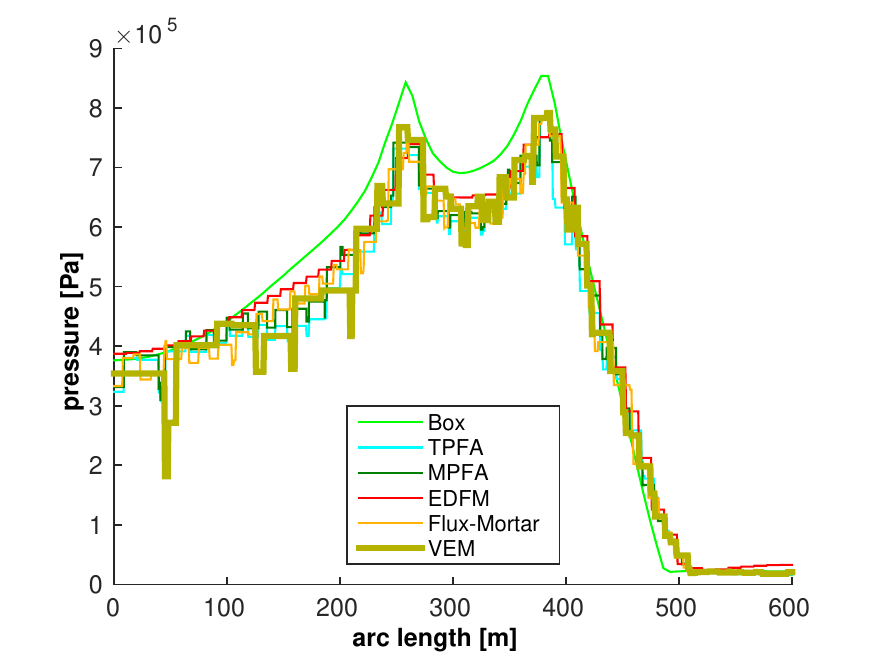}
    \caption{Benchmark 4.
    On the left: pressure solution computed with VEM.
    The pressure ranges in $[0, 1013250]$ Pa. On the centre:
    pressure along horizontal line at $y=500$m. On the
    right: pressure along vertical line at $x=625$m.}%
    \label{fig:sotra_over_line}
\end{figure*}

The computed solution is represented in Figure \ref{fig:sotra_over_line} on the
left, which matches the solution of the others method considered in
\cite{Flemisch2016a}.  Because of the complexity of the network, conforming and
non-matching methods may pose a constraint to the grid generation in particular
for close fractures.  However, the method presented allows general grid cells
lighten the computational cost for the simulation. We refer to Figure
\ref{fig:sotra_grid}, which represents almost intersecting fractures.

As shown in Figure \ref{fig:sotra_grid}, some of the
grid cells are non star-shaped or even contains cuts. For a more detailed
discussion refer to \cite{Fumagalli2017b}.
Finally, to validate in more detail the computed solution we present two plots
over line and compare them with the solutions obtained in the benchmark study,
see Figure \ref{fig:sotra_over_line}.
The curves for the current method are in good agreement with the others.
Small oscillations are related to mesh effects.


\subsection{Passive scalar transport}\label{example_2}

In this part we consider both mixed-dimensional models \eqref{eq:darcy} and
\eqref{eq:transport} to simulate a passive scalar transport. We consider the
geometry presented in Section \ref{subsub:sotra} and compare the solution
obtained with both grids in Figure \ref{fig:sotra_grid} for permeable and
blocking fractures.  For the reference grid the total number of degrees of
freedom is 35578.  The aim of this test is to validate the quality of the Darcy
velocity on the passive scalar in presence of the grid coarsening. In the
following fracture aperture is constant and equal to $10^{-2}$m, a pressure
gradient from the right to left boundary of the domain of $3\cdot 10^7$Pa, a
final simulation time of $40$ years.

\subsubsection{Permeable fractures}

\begin{figure*}[p]
    \centering
    \includegraphics[width=0.33\textwidth]{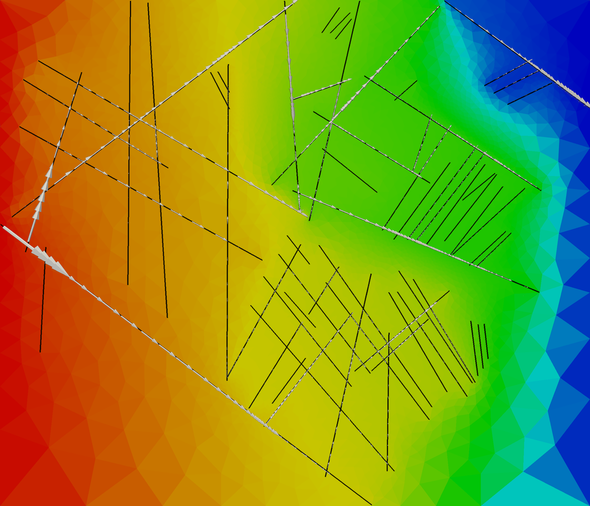}%
    \hspace*{0.006\textwidth}%
    \includegraphics[width=0.33\textwidth]{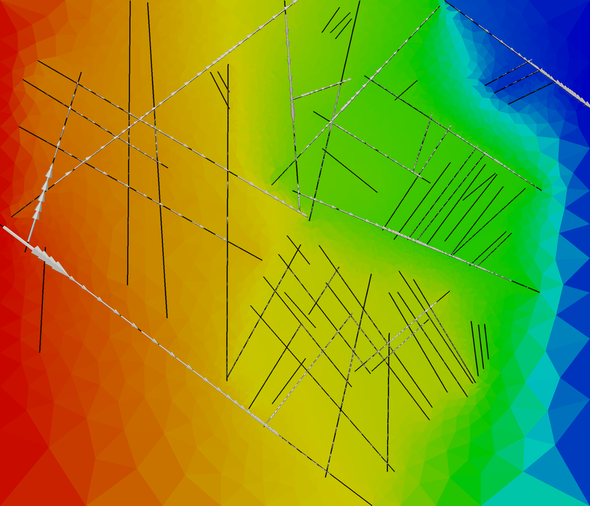}\\
    \includegraphics[width=0.33\textwidth]{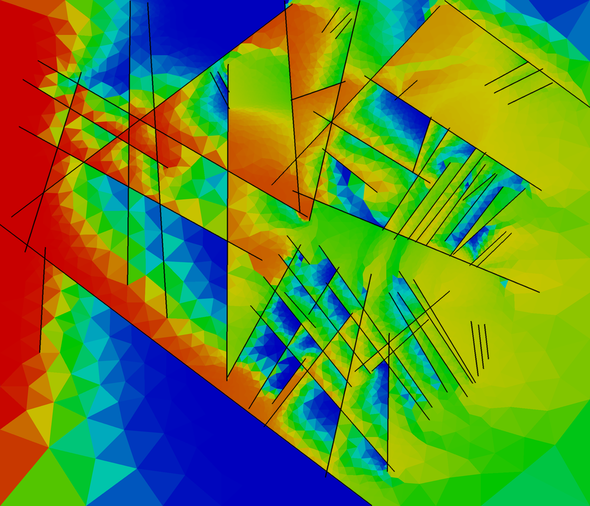}%
    \hspace*{0.006\textwidth}%
    \includegraphics[width=0.33\textwidth]{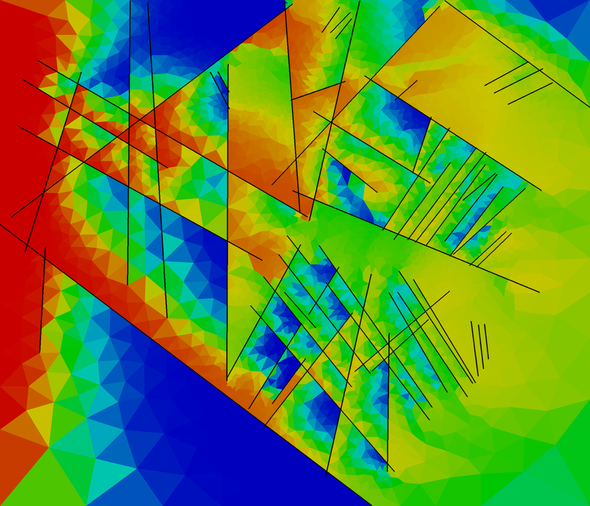}%
    \caption{Permeable fractures.
    On the top: reference and coarse solutions for pressure and
    velocity, as arrows. On the bottom: reference and coarse solutions for
    concentration of the passive scalar.}%
    \label{fig:sotra_transport_permeable}

    \vspace*{\floatsep}

    \centering
    \includegraphics{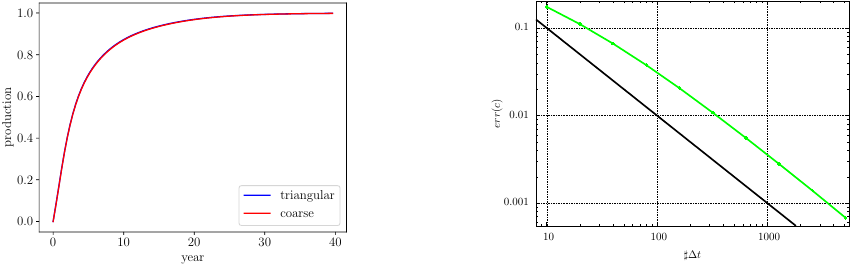}
    \caption{Permeable fractures. On the left: comparison of passive scalar
    production at the outflow between the reference triangular grid and the
    coarse grid. On the right: temporal error decay with reference
    $\mathcal{O}(\Delta t)$ in black.}%
    \label{fig:sotra_transport_permeable_production}
\end{figure*}

We consider high permeable fractures with tangential permeability of $5 \cdot
10^{-6}$ m$^2$ and normal permeability of $2.5\cdot 10^{-9}$ m$^2$. Matrix
permeability is set to $2.5 \cdot 10^{-11}$ m$^2$. Figure
\ref{fig:sotra_transport_permeable} compares the solutions obtained on the
reference triangular grid and on the coarse grid.
Moreover Figure \ref{fig:sotra_transport_permeable_production} on the left
presents a comparison of the passive scalar production. We notice the good
agreement in both the pressure and concentration fields as well as in the
production. We can conclude that in this case the grid coarsening is not
affecting the quality of the computed solutions.

Finally in Figure \ref{fig:sotra_transport_permeable_production} on the right
the temporal error decay is reported. The spatial discretization is fixed and we
consider a sequence of simulation with $(10, 20, 40, 80, 160, 320, 640, 1280,
2560, 5120)$ time steps each.  The error is computed as the L$^2$ difference
from a reference solution obtained with $10^{5}$ time steps after 10 years.  A
unitary error decay is achieved, coherent with the numerical scheme considered.

\subsubsection{Blocking fractures}

\begin{figure*}[p]
    \centering
    \includegraphics[width=0.33\textwidth]{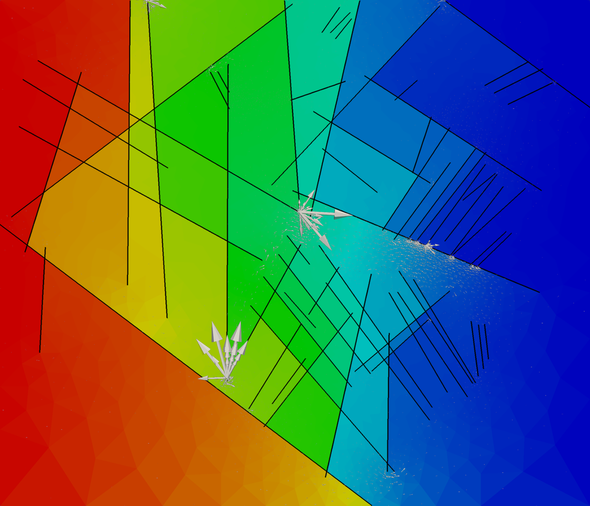}%
    \hspace*{0.006\textwidth}%
    \includegraphics[width=0.33\textwidth]{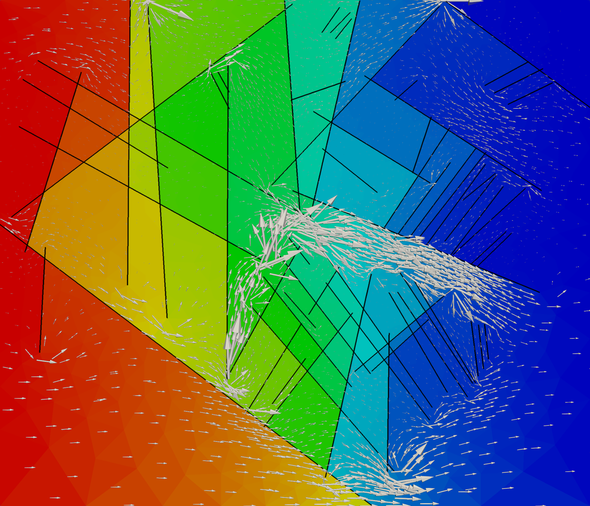}\\
    \includegraphics[width=0.33\textwidth]{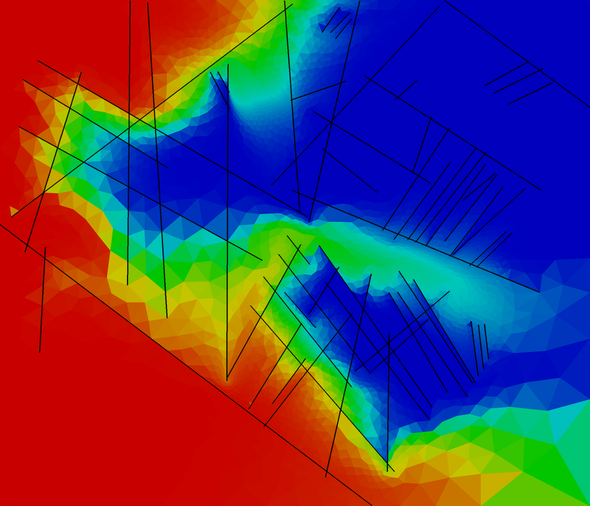}%
    \hspace*{0.006\textwidth}%
    \includegraphics[width=0.33\textwidth]{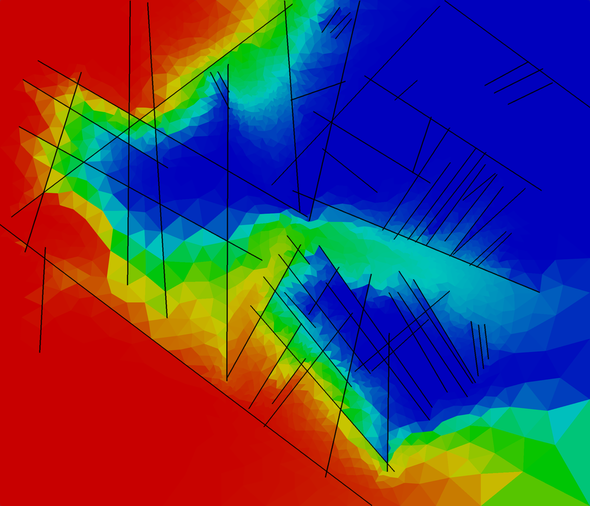}%
    \caption{Blocking fractures.
    On the top: reference and coarse solutions for pressure and
    velocity, as arrows. On the bottom: reference and coarse solutions for
    concentration of the passive scalar.}%
    \label{fig:sotra_transport_blocking}%

    \vspace*{\floatsep}

    \centering
    \includegraphics{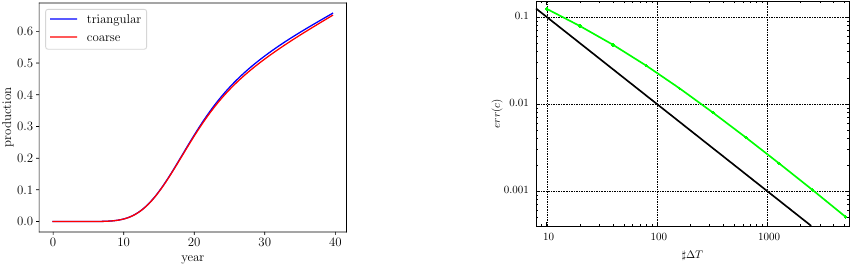}
    \caption{Blocking fractures. On the left: comparison of passive scalar
    production at the outflow between the reference triangular grid and the
    coarse grid. On the right: temporal error decay with reference
    $\mathcal{O}(\Delta t)$ in black.}%
    \label{fig:sotra_transport_blocking_production}
\end{figure*}

We consider low permeable fractures with tangential and normal permeability of $7.5 \cdot
10^{-16}$ m$^2$. Matrix permeability is set to $7.5 \cdot 10^{-11}$ m$^2$. Figure
\ref{fig:sotra_transport_blocking} compares the solutions obtained on the
reference triangular grid and on the coarse grid. In this case we notice a pick
of velocity in the reference case due to small elements close to almost
intersecting fractures. This may affect an explicit in time solver.
Figure \ref{fig:sotra_transport_blocking_production} on the left
presents a comparison of the passive scalar production. We notice the good
agreement in both the pressure and concentration fields as well as in the
production. We can conclude that also in this case the grid coarsening is not
affecting the quality of the computed solutions.

Again, the temporal error decay is reported, see Figure
\ref{fig:sotra_transport_blocking_production}, right. The spatial discretization
is fixed and we require $(10, 20, 40, 80, 160, 320, 640, 1280, 2560, 5120)$ time
steps, respectively. A unitary error decay is achieved, consistent with the
numerical scheme considered.


\subsection{Passive scalar transport on a realistic 3d-network}\label{example_3}

\begin{figure*}[p]
    \centering
    \includegraphics[width=0.145\textwidth]{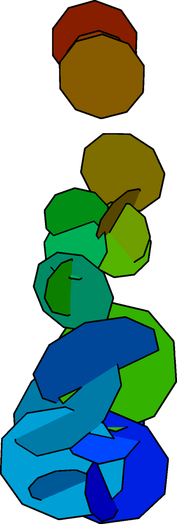}%
    \hspace*{0.05\textwidth}%
    \includegraphics[width=0.145\textwidth]{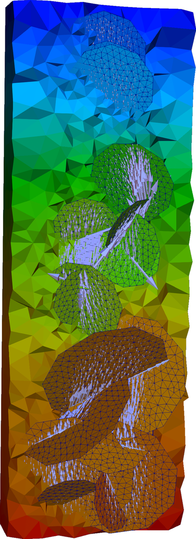}%
    \hspace*{0.05\textwidth}%
    \includegraphics[width=0.145\textwidth]{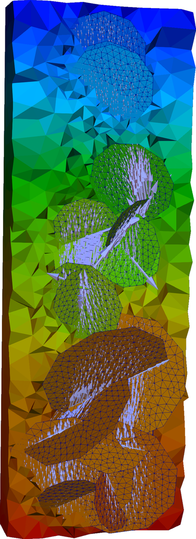}%
    \caption{On the left representation of the 20 fractures coloured by their
        identification number. Pressure and velocity for both the reference and
        coarse grids. The pressure is scaled between 0 and $4.8 \cdot 10^7$.}%
    \label{fig:soultz_pressure}%

    \vspace*{\floatsep}

    \includegraphics[width=0.145\textwidth]{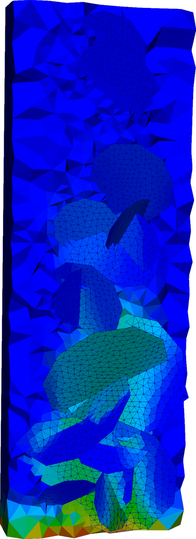}%
    \includegraphics[width=0.145\textwidth]{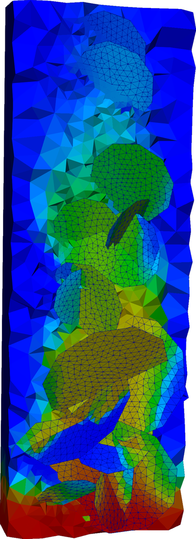}%
    \includegraphics[width=0.145\textwidth]{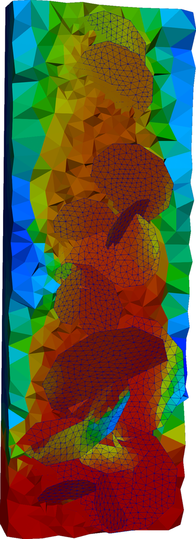}%
    \hspace*{0.05\textwidth}%
    \includegraphics[width=0.145\textwidth]{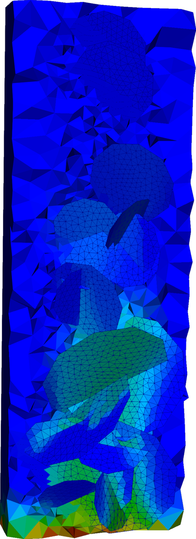}%
    \includegraphics[width=0.145\textwidth]{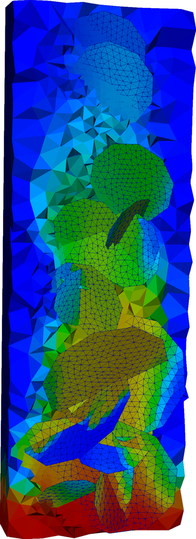}%
    \includegraphics[width=0.145\textwidth]{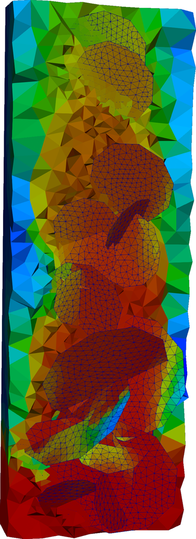}
    \caption{On the left three time steps  on the reference grid for the
    concentration. On the right, at the same steps, the concentration computed
    on the coarse grid.}%
    \label{fig:soultz_transport}%

    \vspace*{\floatsep}

    \includegraphics[width=0.45\textwidth]{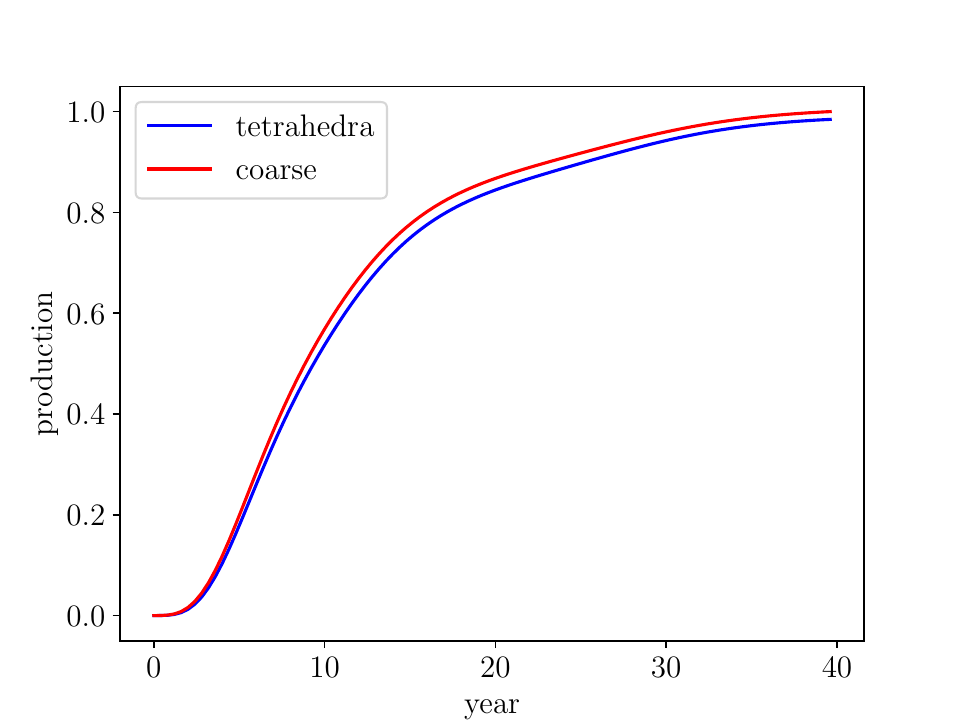}%
    \hspace*{0.05\textwidth}%
    \includegraphics[width=0.45\textwidth]{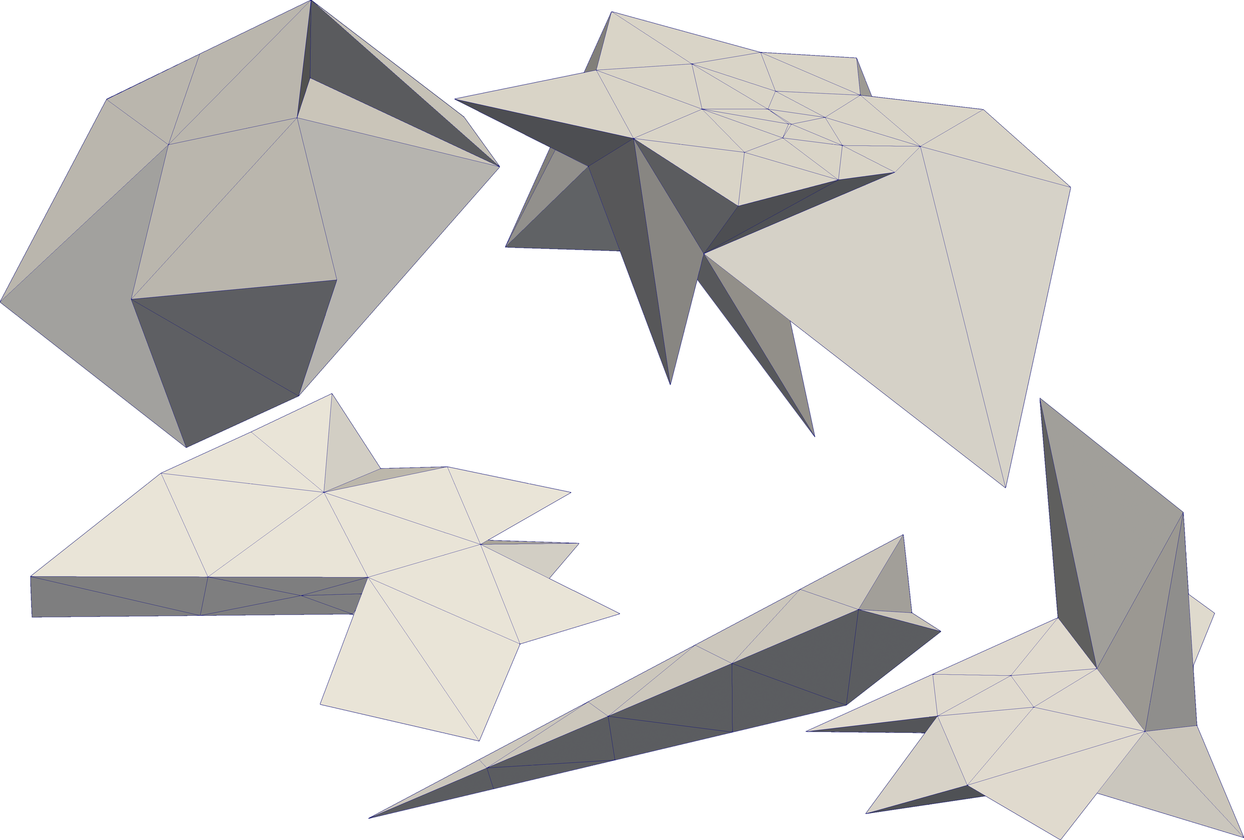}%
    \caption{On the left the production comparison between the reference and the
    coarse solution. On the right an example of five coarse cells.}%
    \label{fig:soultz_production}%
\end{figure*}

In this example we consider a realistic geometry for a geothermal system. We
study a partial reconstruction of fractures from test site at
Soultz-sous-For\^{e}ts in France, for more details see \cite{Sausse2010}. The
network is composed by 20 fractures represented as polygons with 10 edges each.
The fracture intersections result in 33 1d objects and 4 0d objects. In this
case the full model is needed to accurately simulate the fluid flow and
transport on the domain. The fracture geometry is represented in Figure
\ref{fig:soultz_pressure}.
We assume rock matrix permeability as $7.5 \cdot
10^{-10}$m$^2$ and fracture permeability, in both the normal and tangential
direction, equal to $5 \cdot 10^{-5}$m$^2$. The fracture aperture is $\epsilon^2
= 10^{-2}$m
and, for the lower dimensional objects, we consider their ``aperture'' as the
square and the cube of the fracture aperture, respectively for the 1d and 0d
objects. We impose a pressure from the top to bottom of the domain and
no flux boundary conditions on the other sides. The reference grid is composed
by 44331 tetrahedra, 6197 triangles for the 2d grids, 151 segments for the 1d
grids, and 4 point-cells for the 0d grids. After the coarsening algorithm the
resulting grid is composed by 16108 polyhedra, and for the lower dimensional
objects the grids are untouched. See Figure \ref{fig:soultz_production} on the
right for an example of coarse cells.
The transport simulation runs for 40 years and,
for the implicit Euler scheme, we consider 100 time steps.

The objective is to study the robustness of the virtual elements associated with
the coarsening strategy and detect if the coarse model gives accurately enough
results. The pressure solution and the concentration of the scalar passive are
depicted in Figures \ref{fig:soultz_pressure} and \ref{fig:soultz_transport},
respectively, for both the reference and coarse case. Both the pressure and
concentration profiles for the two grids are in good agreement.  To analyse the
macroscopic behaviour of the resulting solution a production curve comparison is
reported in Figure \ref{fig:soultz_production}, showing  a small discrepancy
from the reference and the coarse production.

We can conclude that also in this case
the adopted strategy is effective and can be applied to complex system of
fractures, ensuring a good compromise between high accuracy and computational
effort.

\section{Conclusions} \label{sec:conclusions}

In the paper we presented two classes of mixed-dimensional problems able to
describe a single-phase flow and a transport of a scalar passive in fractured
porous media.  The latter is transported by the Darcy velocity computed by the
former model.  The considered mixed-dimensional Darcy problem is able to
represent channel and barrier behaviours of fractures as well as their
intersections.  The models represent thus a comprehensive description of this
phenomena, which are key ingredients for the description of several energy
application problems. The numerical scheme considered for the discretization of
the Darcy equation are able to handle, supported by theoretical findings, grid
cells of arbitrary geometry becoming a strong advantage when dealing with high
and complex fractured porous media.  Numerical results have shown the capability
to apply this strategy obtaining accurate outcomes with a reasonable
computational cost.  The transport model can be viewed as a first attempt to
introduce a full model for the description of heat exchange in a porous media,
which will be a part of future investigations.




\section*{Acknowledgment}

    We acknowledge financial support for the ANIGMA project from the Research
    Council of Norway (project no. 244129/E20) through the ENERGIX program.
    The authors warmly thank the others components of the
    ANIGMA project team:
    Eivind Bastesen,
    Inga Berre,
    Simon John Buckley,
    Casey Nixon,
    David Peacock,
    Atle Rotevatn,
    P{\aa}l N{\ae}verlid S{\ae}vik,
    Luisa F. Zuluaga.

    The authors wish to thank also:
    Runar Berge,
    Wietse Boon,
    Ivar Stefansson,
    Eren U{\c c}ar.





\end{document}